\documentclass[lfeqn,12pt]{article}
\usepackage{amssymb,bbm,pifont,mathtools,mathabx,mathrsfs,%fdsymbol,
}
\usepackage{minitoc}
 
\usepackage[nottoc]{tocbibind}

\usepackage{comment}

\usepackage{multirow}
\usepackage[multiple]{footmisc}
\usepackage{xcolor}
\usepackage[colorlinks=true]{hyperref}
\newcommand{\giu}{{\medskip\noindent}}

\newcommand{\Giu}{{\bigskip\noindent}}
\newcommand{\nl}{{\smallskip\noindent}}
\newcommand{\noi}{{\noindent}}

\newtheorem{theorem}{Theorem}
\newtheorem{definition}[theorem]{Definition}
\newtheorem{proposition}[theorem]{Proposition}
\newtheorem{lemma}[theorem]{Lemma}
\newtheorem{remark}[theorem]{Remark}
\newtheorem{sublemma}[theorem]{Sublemma}
\newtheorem{corollary}[theorem]{Corollary}
\newtheorem{assumption}[theorem]{Assumption}
\newtheorem{notationalrem}[theorem]{Notational Remark}

\newtheorem{tools}[subsection]{$\negsp\negsp$}
\newcommand\asm[1]{ \begin{assumption}\label{#1} }
\newcommand\easm{ \end{assumption} }
\newcommand\dfn[1]{ \begin{definition}\label{#1} }
\newcommand\dfntwo[2]{ \begin{definition}[#2]\label{#1} }
\newcommand\edfn{ \end{definition} }
\newcommand\rem[1]{ \begin{remark}\label{#1} \small \rm}
\newcommand\remtwo[2]{ \begin{remark}[#2]\label{#1} \rm}
\newcommand\erem{ \end{remark} }
\newcommand\thm[1]{ \begin{theorem}\label{#1}}
\newcommand\thmtwo[2]{ \begin{theorem}[#2]\label{#1}}
\newcommand\ethm{ \end{theorem} }
\newcommand\pro[1]{ \begin{proposition}\label{#1}}       
\newcommand\protwo[2]{ \begin{proposition}[#2]\label{#1}}
\newcommand\epro{ \end{proposition} }
\newcommand\lem[1]{ \begin{lemma}\label{#1}}
\newcommand\lemtwo[2]{ \begin{lemma}[#2]\label{#1}}
\newcommand\elem{ \end{lemma} }
\newcommand\sublem[1]{ \begin{sublemma}\label{#1}}
\newcommand\sublemtwo[2]{ \begin{sublemma}[#2]\label{#1}}
\newcommand\esublem{ \end{sublemma} }
\newcommand\cor[1]{ \begin{corollary}\label{#1}}
\newcommand\cortwo[2]{ \begin{corollary}[#2]\label{#1}}
\newcommand\ecor{ \end{corollary} }
\newcommand\notrem[1]{ \begin{notationalrem}\label{#1} \sl}
\newcommand\enotrem{ \end{notationalrem} }

\newcommand\smalpara{{ \mathfrak{p}}}
\newcommand\average[1]{{ \left\langle #1 \right\rangle}}

\newcommand\equ[1]{{\rm (\ref{#1})}}
\newcommand\beq[1]{ \begin{equation}\label{#1} }
\newcommand{\eeq}{ \end{equation} }

\newcommand\beqa[1]{ \begin{eqnarray} \label{#1}}
\newcommand{\eeqa}{ \end{eqnarray} }
\newcommand{\beqano}{ \begin{eqnarray*} }
\newcommand{\eeqano}{ \end{eqnarray*} }

\newcommand{\proof}{\par\medskip\noindent{\bf Proof\ }}

\newcommand{\ie}{{\it i.e.\  }}

\newcommand{\dst}{\displaystyle}
\newcommand{\qed}{\hskip.5truecm
            \vrule width 1.7truemm height 3.5truemm depth 0.truemm
            \par\Giu}
\newcommand{\qedeq}{\hskip.5truecm
            \vrule width 1.7truemm height 3.5truemm depth 0.truemm}

\newcommand\ovl[1]{ \overline {#1} }

\newcommand\su[1]{ \frac{1}{ {#1}} }
\newcommand\minfoc{ {\, \rm minfoc\, }}
\newcommand{\adj}{ {\rm \, Adj \, }}
\newcommand{\tr}{ {\rm \, tr \, }}

\DeclareMathOperator\dist{dist}

\newcommand\supp{ {\, \rm supp\, }}

\DeclareMathOperator{\meas}{meas}
\DeclareMathOperator{\Lip}{Lip}
\DeclareMathOperator{\diam}{diam}

\newcommand{\diag}{{ \, \rm diag \, }}

\newcommand{\io}{{\infty }}

\newcommand{\dpr}{ {\partial}   }

\newcommand\eqby[1]{\stackrel{\equ{#1}}{=}}
\newcommand\leby[1]{\stackrel{\equ{#1}}{\le}}
\newcommand\ltby[1]{\stackrel{\equ{#1}}{<}}

\newcommand\gtby[1]{\stackrel{\equ{#1}}{>}}
\renewcommand{\Im}{{\rm \, Im\,}}
\newcommand{\inter}{{\rm \, int\,}}

\newcommand{\negsp}{\hspace{-.04truecm}}
\newcommand\ex{\, e}%{\mathbbmss e}}

\renewcommand{\a }{ {\alpha}   }
\renewcommand{\b}{ {\beta}   }
\newcommand{\g}{ {\gamma}   }

\renewcommand{\d}{ {\delta}   }
\newcommand{\D}{ {\Delta}   }
\newcommand{\vae }{ {\varepsilon}   }

\newcommand{\torsion}{\theta}%{\theta}
\renewcommand{\th }{ {\theta}   }

\renewcommand{\k}{\ell}
\renewcommand{\l}{ {\lambda}   }

\newcommand{\e}{ {\varepsilon}   }
\newcommand{\m}{ {\mu}   }
\newcommand{\n}{ {\nu}   }

\newcommand{\p}{ {\pi}   }

\renewcommand{\r}{ {\rho}   }
\newcommand{\s}{ {\sigma}   }

\renewcommand{\t}{ {\tau}   }

\renewcommand{\o}{ {\omega}   }

\newcommand{\torus}{ {\mathbb{ T}}   }
\renewcommand{\natural}{ {\mathbb{ N}}   }
\newcommand{\real}{ {\mathbb{ R}}   }
\newcommand{\integer}{ {\mathbb{ Z}}   }
\newcommand{\complex}{ {\mathbb { C}}   }

\newcommand{\tn}{ {\torus^d} }

\newcommand{\rn}{ {\real^d}   }

\newcommand{\cn}{ {\complex^d }   }
\newcommand{\zn}{ {\integer^d }   }

\newcommand\by{{ \bar y}}

\newcommand{\cH}{ {\cal H} }

\newcommand{\cR}{ {\cal R} }

\newcommand{\cT}{ {\cal T} }

\font\teneufm=eufm10
\font\seveneufm=eufm7
\font\fiveeufm=eufm5
\newfam\eufmfam
\textfont\eufmfam=\teneufm
\scriptfont\eufmfam=\seveneufm
\scriptscriptfont\eufmfam=\fiveeufm

\newcommand\appA[1]{\section{#1}\label{app:A}
\renewcommand{\theequation}{A.\arabic{equation}}
           \setcounter{equation}{0}
\renewcommand{\thetheorem}{A.\arabic{theorem}}
           \setcounter{theorem}{0}
                  }
\newcommand\appB[1]{\section{#1}\label{app:B}
\renewcommand{\theequation}{B.\arabic{equation}}
           \setcounter{equation}{0}
\renewcommand{\thetheorem}{B.\arabic{theorem}}
           \setcounter{theorem}{0}           
           }

\def\bks{{\backslash}}

\def\uno{{\mathbbm 1}}

\newcommand{\wt}{\widetilde}

\newcommand\rball{{\rm B}}
\newcommand\cball{{\,\mathbb B}}

\newcommand{\id}{{\mathtt {id}}}

\newcommand{\actiondom}{{\mathscr D}}
\newcommand{\hactiondom}{{\widehat\actiondom}}

\newcommand\Dioph{{\rm Dioph}}
\newcommand\Do{{\mathcal D_{\!{}_0}}}
\newcommand\DD{{\mathcal D}}
\newcommand\Can{{\actiondom^*}}
\newcommand\Canstar{{\actiondom_*}}

\newcommand\ro{\mathtt{R}}
\newcommand\rr{\mathsf{r}}

\newcommand\KK{{\mathscr K}}

\renewcommand\subset{\subseteq}

\newcommand\cme{\b}

\newcommand\rl{\rho}

\begin{document}

%%%%%%%%%%%%%%%%%%%%%%%%%%%
%%%%%%%   title page

\date{\small \today}

\title{{\bf  V.I. Arnold's ``Global'' KAM Theorem and geometric measure estimates}\\
%\thanks{{\bf Acknowledgments.} We would like to thank Prof. L. Chierchia and Prof. L. Biasco for helping a lot and their patience.}\\
%{\small(Preliminary version)}
}
\author{
L. Chierchia\footnote{Luigi Chierchia: Dipartimento di Matematica e Fisica,  Universit\`a  ``Roma Tre", Largo San Leonardo Murialdo~1, I-00146 Roma (Italy);
\emph{luigi@mat.uniroma3.it}}\ , 
\
C. E. Koudjinan\footnote{Comlan Edmond Koudjinan: Institute of Science and Technology Austria (IST Austria), Am Campus~1, 3400 Klosterneuburg,
Austria; \emph{edmond.koudjinan@ist.ac.at}
} 
}
\maketitle
\begin{abstract}

\noindent
This paper continues the discussion started in \cite{chierchia2019ArnoldKAM} concerning Arnold's legacy on classical KAM theory and (some of) its modern developments. We prove a detailed and explicit `global' Arnold's KAM Theorem, which yields, in particular,  the Whitney conjugacy of a non--degenerate, real--analytic, nearly--integrable Hamiltonian system to an integrable system on a closed, nowhere dense, positive measure subset of the phase space. Detailed measure estimates on the Kolmogorov's set are provided in the case the phase space is: (A) a uniform neighbourhood of an arbitrary (bounded) set times the $d$--torus and
(B) a domain with $C^2$ boundary times the $d$--torus. All constants are explicitly given. 

\end{abstract}

\giu
\section{Introduction}
\begin{itemize}

\item[\bf a.] 
In \cite{chierchia2019ArnoldKAM}, we revised Arnold's original analytic `KAM scheme' \cite{ARV63} and showed, in particular, how to implement it so as to get the optimal relation between the size of the perturbation $\vae$ and the Diophantine constant  
$\a$ associated to a persistent integrable torus (for generalities, we refer to the Introduction in \cite{chierchia2019ArnoldKAM}).

\nl
In the present paper we show how Arnold's `pointwise theorem'  (Theorem A in \cite{chierchia2019ArnoldKAM}) leads, naturally, to a `global theorem', unifying and improving various previous versions of  such a result: compare, in particular, with
\cite{neishtadt1981estimates}, \cite{poschel1982integrability}, \cite{chierchia1982smooth}, \cite{JP}. The term `global' refers, here, to the simultaneous  (and `smooth') construction, in phase space,  of all persistent KAM tori having a prefixed Diophantine constant. The main theorem  (Theorem~\ref{arnoldtheorem} below) is formulated in terms of  a (Whitney) symplectic transformation conjugating a given (Kolmogorov non--degenerate) analytic, nearly--integrable Hamiltonian system  to a Hamiltonian system integrable on a closed, nowhere dense set\footnote{Indeed, closed sets of uniform Diophantine numbers may have, in general, isolated points; compare \cite{Argentieri2020}.}.  All constants involved in Theorem~\ref{arnoldtheorem} are explicitly computed, and, in particular, the optimal relation between $\vae$ and $\a$ is retained.

\item[\bf b.] An immediate corollary of `Arnold's global theorem' is that   {\sl measure estimates} of the (complement of the) Kolmogorov's set (i.e., the set of all persistent integrable tori of a nearly--integrable Hamiltonian system) become essentially trivial (since symplectic transformations preserve Liouville measure on phase space). The problem of finding {\sl explicit measure estimates} of the Kolmogorov's set in terms of the structure of the phase space  is, therefore, reduced to a purely geometrical problem. In particular, as in \cite{biasco2018explicit}, we are interested in analyzing how such measure estimates  depend upon general geometric properties of the action domain, an issue which is particularly relevant in developing KAM theory for {\sl secondary} tori (i.e., those invariant Lagrangian tori which arise because of the perturbation and are not a continuation of integrable tori); compare
\cite{Biasco2015}, \cite{biasco2017secondary},  \cite{biasco2020nonlinearity}. 

\nl
In this paper, we shall discuss detailed measure estimates in two different cases, namely: 

\nl
(A) ({\sl General case}) The Hamiltonian is `uniformly real--analytic'  on  $\DD\times \torus^d$, with  action domain $\DD\subset \real^d$ being  a {\sl completely arbitrary bounded set} and  the unperturbed frequency map is  a local diffeomorphism;
`uniformly analytic' means that the Hamiltonian is real--analytic on the union of complex balls with centers in $\DD$ and {\sl fixed} radius $\ro>0$. In this case the phase space will be $\actiondom\times\torus^d$, where $\actiondom$ is a suitable 
(`minimal') open cover of $\DD$. 
This set--up is similar to that considered in \cite{biasco2018explicit}. 

\nl
(B) ({\sl Smooth case}) The Hamiltonian is real--analytic on a phase space $\actiondom\times\torus^d$ with $\actiondom$  being a bounded, connected, open set  with $C^2$ boundary and the unperturbed frequency map  is a {\sl global} diffeomorphism on $\actiondom$.

\item[\bf c.]
Let us briefly describe the type of measure estimates we get. 

\nl
Case (A): As usual in classical KAM theory, we consider real--analytic Hamiltonians
$$H: (y,x)\in \actiondom\times\torus^d \mapsto 
H(y,x)\coloneqq K(y)+\vae P(y,x)\in\real\,,\eqno{(*)}$$ 
where $(y,x)\in\real^d\times\torus^d$ are standard action--angle variables (i.e., the phase space is endowed with 
the standard symplectic form $dy\wedge dx$), $\vae$ a small parameter,  and  $H$ is   real--analytic on the union of $\ro$--balls with centers in some bounded set $\DD\subset \real^d$, while 
$\actiondom$ is suitable neighbourhood of $\DD$ (see below).
The integrable Hamiltonian $K$ is assumed to be Kolmogorov non--degenerate on $\actiondom$ (i.e., the frequency map
$y\in\actiondom\mapsto \o=\partial_y K(y)$ is a real--analytic local diffeomorphism).
 Let us denote by $\KK_\actiondom(\a,\t)$ the set of Lagrangian graphs over $\torus^d$ in  $\actiondom\times\torus^d$, which are  
invariant for the flow governed by $H$ and on which the flow is analytically conjugated to the Kronecker flow $x\in\torus^d\mapsto x+ \o t$, with $\o\in\real^d$ $(\a,\t)$--Diophantine\footnote{I.e., $|\o\cdot k|\ge \a/|k|^\t$, for any $k\in\integer^d\bks\{0\}$.}, for some $\t>d-1$.  Then, there exist positive numbers $C_*$, $\a_*$, $\vae_*$ and $\rr\le \ro/9$,
depending only on $d$, $\t$,  $K$ and $P$ (and explicitly given in Theorem~4 below), such that, if $0< \vae<\vae_*$, then
$$
\meas\Big((\actiondom\times\torus^d)\backslash \KK_\actiondom(\a_* \sqrt\vae,\t)\Big)\le C_* \,   N_\rr^{\rm int}(\DD)\   \sqrt{\vae}\,,
$$
where  $N_\rr^{\rm int}(\DD)$ is  the so--called $\rr$--internal covering number  of $\DD$ 
and $\actiondom$ is a $\ro$--neighbourhood of a minimal $\ro$--internal cover of $\DD$
(compare \S~\ref{sec:general} for precise definitions).

\nl
Case (B): Here $H$ is as above but $\actiondom$ is assumed to be an open, bounded, connected set with $C^2$ boundary;
$H$ is $\ro$--uniformly real--analytic on $\actiondom$ 
and the unperturbed frequency map is assumed to be a global diffeomorphism on $\actiondom$. Let 
$$
\rr\coloneqq \min\{\ro\,,\,   \minfoc(\partial \actiondom)\,,\, 1/\kappa \}/\sqrt{d}\,,
$$
where `$\minfoc$' denotes the so--called minimal focal distance, and $\kappa$ is the maximum modulus of the principal curvatures of $\partial \actiondom$. Then, there exist positive numbers $\bar C_*$, $\a_*$, and $\vae_*$ 
depending only on $d$, $\tau$, $K$ and $P$ (and explicitly given in Theorem~5 below) 
such that, if $0< \vae<\vae_*$, then
$$
\meas\Big((\actiondom\times\torus^d)\backslash \KK_\actiondom(\a_* \sqrt\vae,\t)\Big)\le \bar C_* \,    \max\big\{\sec_{d-1}(\actiondom)\,,\, \cH^{d-1}(\partial \actiondom)\big\}\     \sqrt{\vae}\,,
$$
where $\sec_{d-1}(\actiondom)$ is the measure of the maximal $(d-1)$--dimensional section of $\actiondom$ and $\cH^{d-1}$ denotes the $(d-1)$--dimensional Hausdorff measure 
(compare \S~\ref{sec:sm} for precise definitions).

\newpage
\item[\bf d.] {\bf Remarks} 

\begin{itemize}

\item[(i)] For the optimality of the relation between $\e$ and $\a$ (and the reason for choosing $\a=\a_* \sqrt{\vae}$ in the Kolmogorov's set), see item {\bf d} in the Introduction of    \cite{chierchia2019ArnoldKAM}.

\item[(ii)] Theorem 4 below extends and generalizes the main result (Theorem~1) in \cite{biasco2018explicit}.

\item[(iii)]
In Appendix A (see, in particular, Remark A.5), we correct a small flaw (concerning the choice of some constants)
in \cite{chierchia2019ArnoldKAM}. 

\item[(iv)]
In Remark A.4 (Appendix A) all constants appearing in the proof are explicitly given. 

\end{itemize}

\item[\bf e.] The paper is organized as follows. 

\nl
In \S~2.1 we introduce some of the notation used in the paper and in \S 2.2 we state the `global Arnold's theorem' (Theorem 1). 
The statement of such a theorem is quite detailed; in particular, the introduction of apparently arbitrary sets or paremeters (such as $\DD_0$ or $\r$) allows to make applications in quite different circumstances (such as  cases (A) and (B) mentioned above). On the other hand, the proof of this theorem does not really contain novel ideas and it is based on the schemes in 
\cite{ARV63}, \cite{koudjinan2019quantitative} and \cite{chierchia2019ArnoldKAM}. However, since we put some emphasis in making everything explicit, we felt necessary to outline the proof, detailing, in particular, the choice of the  (many) parameters involved (this is done in Appendix~A). 

\nl
\S~3 is devoted to measure estimates and, in particular, to the statements and proofs of Theorem 4 and 5, which have been briefly explained in item {\bf c} above.

\nl
Finally, Appendix B contains some of the technical tools used in the paper, namely:
\begin{itemize}
\item[]
{\small
B.1 Classical estimates (Cauchy, Fourier)\\
B.2 An Inverse Function Theorem\\
B.3 Internal coverings\\
B.4 Extensions of Lipschitz continuous functions\\
B.5 Lebesgue measure and Lipschitz continuous map\\ 
B.6 Lipeomorphisms ``close'' to identity\\
B.7 Whitney smoothness\\
B.8 Measure of tubular neighbourhoods of hypersurfaces\\
B.9 Kolmogorov non--degenerate normal forms\\
}
\end{itemize}

\end{itemize}

\section{Arnold's Global KAM Theorem}
\subsection{Notations}
\label{parassnot}
\begin{itemize}
\item[\tiny $\bullet$]$\natural\coloneqq\{1,2,3,\cdots\}$ and $\natural_0\coloneqq\{0,1,2,3,\cdots\}$.

\item[\tiny $\bullet$] For $d\in\natural $ and  $x,y\in \complex^d$, we let
$x\cdot y\coloneqq x_1 \bar y_1 +\cdots+x_d \bar y_d$ be the standard inner product (the bar denotes complex conjugate). We denote, respectively, the sup--norm, the 1--norm and the Euclidean norm, by:
$$
|x|\coloneqq \dst\max_{1\le j\le n}|x_j|\,,\qquad  
|x|_1\coloneqq\sum_{j=1}^d |x_j|\,, 
\qquad
|x|_2\coloneqq\sqrt{\sum_{j=1}^d |x_j|^2}\,.
$$

\item[\tiny $\bullet$] $\tn\coloneqq \rn/2\p\zn$ is the  $d$--dimensional (flat) torus.

\item[\tiny $\bullet$] Given $\a>0$, $\t\ge d-1\ge 1$,  we denote by
\beq{dio}\Dioph_\a^\t\coloneqq \big\{\o\in \rn:|\o\cdot k|\geq \frac{\a}{|k|_1^\t}, \ \ \forall\ 0\not=k\in \zn\big\},
\eeq
the set of {\sl $(\a,\t)$--Diophantine vectors} in $\real^d$.

\item[\tiny $\bullet$] For $r,s>0$, $y_0\in\complex^d$, $\emptyset\neq D\subseteq \cn$, we denote:
\beqano
\rball_r(y_0)&\coloneqq  &\left\{y\in \real^d: |y-y_0|<r\right\}\,,\qquad \phantom{A}(y_0\in\real^d)\,,\\
\rball_r(D)&\coloneqq  &\dst\bigcup_{y_0\in D}\rball_{r}(y_0)\,,\phantom{AAAAAAAAAAa}\;(D\subset \real^d)\,,\\
\cball_r(y_0)&\coloneqq  &\left\{y\in \cn: |y-y_0|<r\right\}\,,\\
\cball_{r}(D)&\coloneqq & \dst\bigcup_{y_0\in D}\cball_{r}(y_0)\,,\\
\dst\torus^d_s &\coloneqq  &\left\{x\in \cn: |\Im x|<s\right\}/2\p\zn\,,\\
\cball_{r,s}(y_0)&\coloneqq  & \cball_r(y_0)\times \torus^d_s\,,\\
\cball_{r,s}(D)&\coloneqq  & \cball_r(D)\times \torus^d_s\,;
\eeqano 
we shall also denote, in bold face characters, Euclidean balls:
\beqano
{\bf B}_r(y_0)&\coloneqq  &{\left\{y\in \real^d: |y-y_0|_2<r\right\}\,,\qquad (y_0\in\real^d)\,,}\\
{\bf B}_r(D)&\coloneqq  &\dst\bigcup_{y_0\in D}{\bf B}_{r}(y_0)\,,\phantom{AAAAAAAAAaa}(D\subset \real^d)\,.
\eeqano

\item[\tiny $\bullet$] If  $\uno_d\coloneqq \diag(1)$ is the unit $(d\times d)$ matrix, we denote the standard symplectic matrix by 
$$\mathbb{J}\coloneqq  \begin{pmatrix}0 & -\uno_d\\
\uno_d & 0\end{pmatrix}\,.
$$ 

\item[\tiny $\bullet$]
For $D\subset \rn$, $r\ge 0$ and $s>0$, $\mathcal{B}_{r,s}(D)$ denotes the Banach space of real--analytic functions 
$$f:\cball_r(D)\times \torus_s\to\complex$$ with bounded holomorphic extensions to $ \cball_{r,s}(D)$, with  uniform norm
$$
\|f\|_{r,s}\coloneqq \|f\|_{r,s,D}\coloneqq \dst\sup_{ \cball_{r,s}(D)}|f|<\io\,.
$$
Analogously, $\mathcal{B}_{r}(D)$ denotes the Banach space of real--analytic functions 
$$f:\cball_r(D) \to\complex$$ with bounded holomorphic extensions to $ \cball_{r}(D)$, with 
$$
\|f\|_r\coloneqq \|f\|_{r,D}\coloneqq \dst\sup_{ \cball_{r}(D)}|f|<\io\,.
$$

\item[\tiny $\bullet$]
For a differentiable function $f\colon A\subset  \cn\times\cn \ni (y,x)\mapsto f(y,x)\in \complex$, its gradient/Jacobian is denoted by $\nabla f$ or by $f'$. 

\item[\tiny $\bullet$] We equip 
$\cn\times\cn$   (and its subsets) with the canonical symplectic form 
$$\varpi\coloneqq dy\wedge dx=dy_1\wedge dx_1+\cdots+dy_d\wedge dx_d\ ,
$$
and denote by $\phi_H^t$ the associated Hamiltonian flow governed by the Hamiltonian $H(y,x)$, $y,x\in\complex^d$, i.e., $z(t)\coloneqq \phi_H^t(z)$ is the unique solution of 
$$\dot z= \mathbb{J} \nabla H\,,\qquad  z(0)=z\,.$$

\item[\tiny $\bullet$]
Given a linear operator $L$ from the normed space $(V_{\rm a},\|\cdot\|_{\rm a})$ into the normed space  $(V_{\rm b},\|\cdot\|_{\rm b})$, its ``operator--norm'' is given by
\[\|L\|\coloneqq \sup_{x\in V_{\rm a}\setminus\{0\}}\frac{\|Lx\|_{\rm b}}{\|x\|_{\rm a}},\quad \mbox{so that}\quad \|Lx\|_{\rm b}\le \|L\|\, \|x\|_{\rm a}\quad \mbox{for any}\quad x\in V_{\rm a}.\]
\item[\tiny $\bullet$]
Given $\o\in \rn$, the directional derivative of a $C^1$ function $f$ with respect to $\o$ is given by
 
\[D_\o f\coloneqq \o\cdot  f_x=\dst\sum_{j=1}^d \o_j \dst f_{{x}_j}\,.\]

\item[\tiny $\bullet$]
If  $f$ is a (smooth or analytic)  function on $\torus^d$, its Fourier expansion is given by    \[f=\dst\sum_{k\in \zn}f_k \ex^{ik\cdot x}\,,
\qquad f_k:=\dst\frac{1}{(2\pi)^d}\dst\int_{\tn}f(x) \ex^{-ik\cdot x}\, d x\,,\]
(where, as usual, $\ex\coloneqq \exp(1)$ denotes the Neper number and $i$ the 
imaginary unit). We also set:
\[\average{f}\coloneqq f_0=\dst\frac{1}{(2\pi)^d}\dst\int_{\tn}f(x)\, d x\,,\qquad  T_N f:=\dst\sum_{|k|_1\leq N}f_k \ex^{ik\cdot x},\, N>0. 
\]
\item[\tiny $\bullet$]
For a function $f\colon (\mathscr M_1,{\rm d}_1)\to (\mathscr M_2,{\rm d}_2)$, where $(\mathscr M_j,{\rm d}_j)$, $j=1,2$ are metric spaces, we denote 
$$
{\Lip_{\mathscr M_1}(f)}\coloneqq\|f\|_{L,\mathscr M_1}\coloneqq \sup_{x\neq x'\in \mathscr M_1}\frac{{\rm d}_2(f(x),f(x'))}{{\rm d}_1(x,x')}\le \infty,
$$
 and $f$ is said Lipschitz continuous on $\mathscr M_1$  if $\Lip_{\mathscr M_1}(f)<\infty$.\\
  If $\mathscr M_1=\rn$, we usually 
 denote $\Lip_\rn(f)=\Lip(f)$.

 \item[\tiny $\bullet$]
 $C^k_W(D)$ denotes the set of functions which are $C^k$ in the sense of Whitney on the set $D$. A $C^1_W$  map
 $\phi:D\times\tn \to \rn\times \torus^d$, is symplectic if
 the Whitney--gradient $\nabla \phi=(\dpr_{y}\phi,\dpr_x \phi)$ satisfies $(\nabla \phi)\mathbb{J}(\nabla \phi)^T=\mathbb{J}$  on $D\times\tn$. For more details,  see Appendix~\ref{appD}.

 \item[\tiny $\bullet$]
 The $s$--dimensional Hausdorff measure on $\rn$ will be denoted by $\mathcal{H}^s$; in particular  $\mathcal{H}^d$, which coincides with the $d$--dimensional outer Lebesgue measure, will be denoted by `$\meas$'.
\end{itemize}

\subsection{KAM Theorem}
Given an open set $\actiondom\subset \real^d$ and a real--analytic Hamiltonian $H:\actiondom\times\torus^d\to \real$, we say that {\sl $\cT\subset \actiondom\times \torus^d$ is a (primary\footnote{As opposed to {\sl secondary tori} (same definition but removing the graph assumption); for a KAM Theory for secondary tori, see \cite{Biasco2015}. In this paper, we shall only consider primary KAM tori})  Kolmogorov (or `KAM')  torus for  $H$} if $\cT$ is a real--analytic Lagrangian embedded torus  $\cT=\phi(\torus^d)$, which is a graph over $\torus^d$, and such that  
$$\phi_H^t\big( \phi (\th)\big) = \phi(\th+\o t)\,,\qquad \forall\ \th\in \torus^d\,, t\in\real\,,$$
 for a given Diophantine `frequency vector' $\o\in \Dioph_\a^\t$ (for some $\a,\t>0)$.

\thm{arnoldtheorem}
Let $d \ge 2$;   $\ro>0$; $0<s\le 1$;  $\emptyset \neq \DD\subset \rn$; $\vae,\a>0$. Let 
the `integrable Hamiltonian' 
$K\in \mathcal{B}_\ro(\DD)$ be a uniformly (Kolmogorov) non--degenerate   (i.e. $\det K_{yy}\neq 0$ on 
$%\overline
{\cball_\ro(\DD)}$) and let the `perturbation' $P$ belong to 
 $\mathcal{B}_{\ro,s}(\DD)$. Define
\beqa{defpar}
&&
\mathsf{M} \coloneqq \|K_{yy}\|_{\ro,\DD}\;,\quad
\mathsf{L} \coloneqq  \|K_{yy}^{-1}\|_{\ro,\DD}\,,\quad
\mathsf{P} \coloneqq \|P\|_{\ro,s,\DD}\;, \quad\torsion \coloneqq \mathsf{M} \mathsf{L}\,, 
%&& 
\quad
\epsilon\coloneqq \vae\, \frac{ \mathsf{M}  \mathsf{P}}{\a^2}\,.
\eeqa
Choose $0<\rl<\rr\le \ro$, $\Do\subset \DD$, $\t\ge d-1$; define the following `action domains':
\beq{actiondomains}
\actiondom\coloneqq \rball_\rr(\Do)\,,\qquad 
\hactiondom\coloneqq  \rball_{\rr-\rl}(\Do)\,,\qquad
\Can\coloneqq \big\{y\in \hactiondom:\ K_y(y) \in \Dioph_\a^\t\big\}\,,
\eeq
and consider the `nearly--integrable',  non--degenerate Hamiltonian  given by
$$H: (y,x)\in \actiondom\times\torus^d \mapsto 
H(y,x)\coloneqq K(y)+\vae P(y,x)\in\real\,;$$ 
the `phase space' $\actiondom\times\torus^d$ being endowed with the standard symplectic form $\varpi$. Fix $0<s_*<s$.

\nl
There exist constants $c_*, c_0, c_1,  c_2, c_3, c_4>1$,  depending only on $d$ and $\t$, such that, if  
\beq{smcEAr0v2}
\a\le c_0\;\frac{\rl }{\mathsf{L}}\,;
\phantom{AAAAAA}  \epsilon\le \epsilon_*\coloneqq\frac{(s-s_*)^a}{c_*\;\torsion^6}\,,
\eeq 
with {$a\coloneqq 7\n+4d+2$} and $\n\coloneqq \t+1$, 
then, the following statements holds. \\
There exists a nowhere dense set 
$\Canstar\subset \rball_{\rr-\frac\rl2}(\Do)\subset {\actiondom}$,
a lipeomorphism 
$$Y^*\colon \Can\overset{onto}{\longrightarrow}\Canstar\,,$$ 
a function
$K_*\in C_W^\infty(\Canstar)$
and  a $C_W^\infty $--symplectic  transformation 
\beq{KolmogorovSet}
\phi_*\coloneqq \id+(v_*,u_*)\colon \Canstar\times \tn\to \KK\coloneqq \phi_*(\Canstar\times \tn)\subset \actiondom\times \tn\,,
\eeq
real--analytic in $x\in\torus^d_{s_*}$,
 such that\footnote{$y_*$--derivatives are  Whitney--derivatives.} 
\begin{align}
\dpr_{y_*}K_*\circ Y^*&=\dpr_{y}K\,,  \phantom{AAAAAa} \mbox{on} \quad \Can\;,\label{conjCaneq00v2}\\
\dpr^\b_{y_*}(H\circ \phi_*)(y_*,x)&=\dpr^\b_{y_*}K_*(y_*),\qquad \forall\;(y_*,x)\in \Canstar\times\tn,\quad \forall\; \b\in\natural_0^d \label{conjCaneq0v2}\,.
\end{align}
Furthermore, the following estimates hold:
\begin{align}
&\|Y^*-\id\|_{\Can}
\le 
{c_1}  
\;(s-s_*)^\n\;\torsion^2\;\frac{ \vae  {\mathsf{P}}}{\a}\,,\label{NormGstrThtv20}\\
& 
{\Lip_\Can(Y^*-\id)}\le 
{c_2} \;\torsion^3\;\;(s-s_*)^{-1}\;\frac{\mathsf{M} {\vae}\mathsf{P}}{\a^2}\;\left(\log\frac{\a^2}{\mathsf{M} {\vae}\mathsf{P}}\right)^\n\le \su{4d}\,,\label{NormGstrThtv21}\\
& \max\left\{\|u_*\|_*,\; 2d\sqrt{2}\,\frac{\mathsf{M}\k^\n}{\a}\|v_*\|_* \right\}\le  
{c_3}
\;\k^{\n}\;\torsion^2\;\frac{\mathsf{M} {\vae}\mathsf{P}}{\a^2},\label{estArnTrExtv20}\\
& \|\dpr_x u_*\|_*\le 
{c_4} 
\;\torsion_0^2\;\k^{\n}\;\frac{\mathsf{M} {\vae}\mathsf{P}}{\a^2}\le \su{4(18d^3+70)\torsion} \;,\label{estArnTrExtv21}
\end{align}
where
$$
\|\cdot\|_*\coloneqq \sup_{\Canstar\times\torus^d_{s_*}}|\cdot|, \qquad 
{\k\coloneqq 8(s-s_*)^{-1}\log\epsilon^{-1}}
\;.
$$
The `Kolmogorov set' $\KK$  defined in \equ{KolmogorovSet} is foliated, as $y_*\in \Canstar$, by Kolmogorov tori
{$\cT_*\coloneqq\phi_*(\{y_*\}\times \tn)$},  which  are Kolmogorov non--degenerate\footnote{For the precise definition, see  Appendix~\ref{app:A} and \ref{app9}.}.
\ethm
The proof of this theorem is based upon Arnold's original KAM scheme, revised and improved in \cite{chierchia2019ArnoldKAM}, where, in particular all constants are computed and optimal smallness conditions concerning the relation between small divisors and smallness of the perturbation are given. Since, essentially, no new ideas are needed, details are deferred to Appendix~\ref{app:A}. 

\nl
However, let us make, here, a few observations. 

\rem{rem:Thm1} (i) The hypotheses on $H$ can be rephrased by saying that $H$ is $\ro$--uniformly real--analytic on $\DD$.
Notice that $\DD$ can be a completely arbitrary subset of $\rn$, but $\actiondom$ and $\hactiondom$ are open sets.
\\
(ii) The introduction of $\Do$ and $\rl$ is made in order to be able to apply the theorem in quite different contexts; compare, e.g.,  next section on measure estimates. 
\\
(iii) Even if $\Do$ is a single point, the theorem guarantees, in general,  a set of positive measure of Kolmogorov tori for $H$, since the set $\Can$ is a set of positive measure, provided  $\t>d-1$ and $\a$ is small enough. Precise measure estimates are one of the  objectives of this paper and will be given in next section.
\\
(iv) The parameter  $\torsion$ defined in \equ{defpar} measures the `torsion' of the unperturbed system and it is always greater or equal than 1; indeed,  for any $y_0\in\actiondom$, denoting $T(y)\coloneqq K_{yy}(y)^{-1}$,  one has:
\beq{Eta0Ge1}
\torsion\coloneqq \mathsf{L}\mathsf{M}\ge \|T(y_0)\|\|K_{yy}(y_0)\|=\|T(y_0)\|\|T(y_0)^{-1}\|\ge 1\ .
\eeq
\\
(v) The constants $c_i$ appearing in the theorem are explicitly given in  Appendix~\ref{app:A}; compare, in particular,  Eq.~\equ{constants}.
\\
(vi) $\Can$ and $\Canstar$ are closed, nowhere dense sets, but may have isolated points\footnote{See \cite{Argentieri2020}.}.
\erem

\section{Measure estimates\label{MeSeS}}
The fact that the Kolmogorov set $\KK$ in Theorem~\ref{arnoldtheorem} is the image of a (Whitney) symplectic map, leads to straightforward measure estimates of its complement:

\thm{thmmeas1} Under the same notations and assumptions of Theorem~\ref{arnoldtheorem}, 
let
\beqa{defsets}
&& \cme\coloneqq\big(1+ 2\Lip_\Can(Y^*-\id)\big)^d\, (2\pi)^d\,,
\nonumber
\\
&& \mathscr T_\rl\coloneqq\rball_{\rr+\rl}(\Do)\setminus\rball_{\rr-\rl}(\Do) \,,\nonumber\\
&&
\cR_\a\coloneqq\big\{y\in \actiondom:K_y(y)\notin \Dioph_\a^\t\big\}\,.
\eeqa
Then, one has 
\beqa{MesChPin}
\meas(\actiondom\times\tn\setminus \KK)&\le&
\cme \,  
\meas\big(\rball_{\frac\rl2 }(\actiondom)\setminus \Can\big)\\
 &\le& \cme\big(
\meas({\mathscr T}_\rl)+\meas(\cR_\a)\big)
\,.\nonumber
\eeqa
\ethm

\proof
By Theorem~\ref{Minty}, we can extend $Y^*-\id$ component--wise to obtain a global Lipschitz continuous function $f\colon \rn\righttoleftarrow$ satisfying $f|_{\Can}=Y^*-\id$ and
\begin{align}
\dst\sup_{\rn}|f|&=\dst\sup_{\Can}|Y^*-\id|\stackrel{\equ{NormGstrThtv20},\equ{smcEAr0v2}}{\le} \frac\rl2\;,
\qquad \Lip_\rn(f)=  \Lip_\Can(Y^*-\id)<\frac{1}{4d}\;.\label{LipgloGstr0K}
\end{align}
\begin{comment},
Hence,
\begin{align}
\|f\|_{\rn}&\overset{def}{=}\dst\sup_{\rn}|f|\nonumber\\
		   &\le \dst\sup_{\rn}|f|_{2}\nonumber\\
	       &\eqby{LipgloGstr0}\dst\sup_{\Can}|Y^*-id|_{2}\nonumber\\
	       &\le \sqrt{d}\dst\sup_{\Can}|Y^*-id|\nonumber\\
	       &\leby{NormGstrThtv2}\sqrt{d}\;r_*\nonumber\\
	       &\ltby{smcEAr0v2}\su{32d}\frac{r_0\s_0}{\torsion_0}\nonumber\\
	       &\le {\d \s_0} \label{LipgloGstr01K}
\end{align}
and
\begin{align}
\|f\|_{L,\rn}&\overset{def}{=}\dst\sup_{\substack{y,y'\in\rn\\y\neq y'}}\frac{|f(y)-f(y')|}{|y-y'|}\nonumber\\
		     &\le \dst\sup_{\substack{y,y'\in\rn\\y\neq y'}}\frac{|f(y)-f(y')|_{2}}{|y-y'|_2/\sqrt{d}}\nonumber\\
		     &\eqby{LipgloGstr1K}\sqrt{d}\dst\sup_{\substack{y,y'\in\Can\\y\neq y'}}\frac{|(Y^*-id)(y)-(Y^*-id)(y')|_{2}}{|y-y'|_2}\nonumber\\
		     &\le \sqrt{d}\dst\sup_{\substack{y,y'\in\Can\\y\neq y'}}\frac{\sqrt{d}|(Y^*-id)(y)-(Y^*-id)(y')|}{|y-y'|}\nonumber\\
		     &=d\|Y^*-\id\|_{L,\Can}\nonumber\\
		     &\leby{LipGstrThtv2} d\frac{\ex\;\s_0^{\n+d}}{\mathsf{C}_6}<\su 2\;.\label{LipgloGstr11K}
\end{align}
\end{comment}
Set $g\coloneqq f+\id$. Then, by Lemma~\ref{LipRang} and \equ{LipgloGstr0K}, one has\footnote{The bar on sets denotes closure.} 
\beq{GsDdeltMT}
\actiondom\subset g\big(\ovl{ B_{\frac\rl2}(\actiondom)}\big)
\;.
\eeq 
Notice also that, by \equ{LipgloGstr0K}
and Lemma~\ref{LipRang}, 
$g$ is a lipeomorphism of $\rn$. Consequently,
\beqano
\meas(\actiondom\times\tn\setminus \KK)&\coloneqq&\meas(\actiondom\times\tn)-\meas\big(\phi_*(\Canstar\times\tn)\big)\\
&=&\meas(\actiondom\times\tn)-\meas(\Canstar\times\tn)\\
&=&(2\pi)^d \big(\meas(\actiondom)-\meas(\Canstar)\big)\\
&\leby{GsDdeltMT}&(2\pi)^d \big( \meas\big(g(\ovl{ B_{\frac\rl2}(\actiondom)})\big)-\meas(\Canstar)\big)\\
 &=&(2\pi)^d \meas\big(g(\ovl{ B_{\frac\rl2}(\actiondom)})\setminus g(\Can)\big)\\
 &=&(2\pi)^d\meas\big(g\big(\ovl{ B_{\frac\rl2}(\actiondom)}\setminus \Can\big)\big)\qquad\qquad \mbox{(because $g$ is injective)}\\
 &\stackrel{\equ{upperbound}}{\le}&(2\pi)^d(\Lip g)^d\meas\big(\ovl{ B_{\frac\rl2}(\actiondom)}\setminus \Can\big)\\
 &\leby{LipgloGstr0K}&(2\pi)^d(1+ 2\Lip(Y^*-\id))^d\meas\big({B}_{\frac\rl2}(\actiondom)\setminus \Can\big)\,.
\eeqano 
Finally, recalling that $\actiondom=\rball_\rr(\Do)$ and \equ{actiondomains}, one sees that
\beqano
\rball_{\frac\rl2 }(\actiondom)\setminus \Can&=&\rball_{\frac\rl2}(\actiondom)\setminus \hactiondom\  \dot \cup\  \hactiondom\setminus\Can\\
&{=}&{\rball_{\rr+\frac\rl2}(\Do)\setminus\rball_{\rr-\rl}(\Do)\  \dot \cup\  \big\{y\in \rball_{\rr-\rl}(\Do):K_y(y)\notin \Dioph_\a^\t\big\}}\\
&\subset & {\mathscr T}_\rl \cup \cR_\a\,,
\eeqano
from which, the second inequality in \equ{MesChPin} follows at once. \qed

\nl
Theorem~\ref{thmmeas1} reduces the problem of estimating the measure of the complement of the Kolmogorov set $\KK$  to the  estimate on the measure of the complement of Diophantine numbers in a given set and to the purely geometrical problem of estimating the measure of the {\sl tubular neighbourhood} ${\mathscr T}_\rl$ of the boundary of $\actiondom=\rball_\rr(\Do)$. 
Therefore, concrete measure estimates will depend upon the structure of the action domain $\actiondom$ and of 
the (unperturbed) frequency map 
\beq{frequencymap}
y\in \actiondom\mapsto \o_0(y)\coloneqq K_y(y)\in\rn\,. 
\eeq

\nl
We shall discuss, in detail, two different cases:

\begin{itemize} 

\item[(A)] ({\sl General case}) $\DD$ is an {\sl arbitrary bounded} set, $H$ uniformly  real--analytic on $\DD\times \torus^d$  and $\o_0$ is a {\sl local} diffeomorphism on $\DD$ (which is always the case, if   the unperturbed Hamiltonian is assumed  to be Kolmogorov non--degenerate). In this case, as phase space we shall consider a  `minimal' (in a suitable sense) open cover of $\DD$ times $\torus^d$. This set--up is analogous to that considered in \cite{biasco2018explicit}.

\item[(B)] ({\sl Smooth case}) $\DD$ is a bounded, connected, open set  with $C^2$ boundary and $\o_0$ si a {\sl global} diffeomorphism on $\DD$. In this case the phase space is just $\actiondom\times\torus^d\coloneqq \DD\times \torus^d$.
\end{itemize}

\subsection{General case}\label{sec:general}
In order to state the result for  case (A), let us give two definitions.

\begin{itemize}

\item[\tiny $\bullet$]
Given a bounded non--empty set $\DD\subset \real^d$, and given $r>0$,  an {\bf $r$--internal covering of $\DD$} is a 
subset $\Do$ of $\DD$ such that 
\beq{covered}
\DD\subset \rball_r(\Do)=\bigcup_{y\in \Do} \rball_r(y)\,;
\eeq
 $N_r^{\rm int}(\DD)$, the {\bf $r$--internal covering number  of $\DD$}, is defined as\footnote{$N_r^{\rm int}(\DD)$ is finite if and only if $\DD$ is bounded. A simple upper bound on $N_r^{\rm int}(\DD)$ for bounded domains $\DD$ is:
$N_r^{\rm int}(\DD)\le ([\diam(\DD)/r]+1)^d$; compare \cite{biasco2018explicit} or Appendix~{B}, \S~\ref{app:covering}.}
\beq{coveringnumber}
N_r^{\rm int}\coloneqq \min \big\{n\in\natural:\ \{y_1,...,y_n\}\ {\rm is\ an}\  r\!-{\!\rm internal\ covering\ of\ } \DD\big\}\,;
\eeq
an $r$--internal cover $\Do$ of $\DD$ with cardinality equal to the $r$--internal covering number will be called {\bf a minimal $r$--internal covering of $\DD$}. 

\item[\tiny $\bullet$]
Given a real--analytic Hamiltonian $H:\actiondom\times\torus^d\to\real$, 
we denote the set of KAM tori for $H$ with frequency vector in $\Dioph_\a^\t$ by
\beq{KAMtori}
\KK_\actiondom(\a,\t)\coloneqq\big\{ \cT\subset \actiondom\times\torus^d|\,\cT\ {\rm is\ a\ KAM\ torus \ for}\  H\ {\rm with\ 
frequency}\ \o\in \Dioph_\a^\t\big\}.
\eeq
\end{itemize}

\thm{thmmeas2} Let $\DD$ be an arbitrary bounded non--empty set in $\real^d$, $\t>d-1\ge1$, $\ro,s>0$. Let 
 $K\in \mathcal{B}_\ro(\DD)$ be a uniformly (Kolmogorov) non--degenerate, $P\in\mathcal{B}_{\ro,s}(\DD)$ and let $\mathsf{M}$, $\mathsf{L}$, $\mathsf{P}$, $\torsion$ as in \equ{defpar}. Let $c_0$ and $c_*$ be as in Theorem~\ref{arnoldtheorem}. 
 Fix $0<s_*<s$, let $\epsilon_*$ be as in \equ{smcEAr0v2} and  define:
 \beqa{def.meas1}
&&\rr\coloneqq \frac{\ro}{1+2d^2\torsion}\,,\qquad \a_*\coloneqq \sqrt{\frac{\mathsf{M}\mathsf{P}}{\epsilon_*}}\,,
\qquad  \vae_*\coloneqq \Big(\frac{c_0 \rr}{\mathsf{L}\a_*}\Big)^2\,,
\nonumber\\
&& \d_0\coloneqq \inf_{\rball_\rr(\DD)}|\det K_{yy}|\,,\quad \torsion_0\coloneqq \max\Big\{ \frac{\mathsf{M}^d}{\d_0}\,,\torsion\Big\}\,,\quad \r\coloneqq \frac{\a_* \mathsf{L}}{c_0}\, 
\sqrt\vae\,.
\eeqa
Let $\Do\subset \DD$ be a minimal  $\rr$--internal covering of $\DD$, $\actiondom\coloneqq \rball_\rr(\Do)$ and  let
$\KK_\actiondom(\a_* \sqrt\vae,\t)$ be as in \equ{KAMtori} with $H=K+\vae P$.
Then, if $0< \vae<\vae_*$, 
one has
\beq{eas1}
\meas\Big((\actiondom\times\torus^d)\backslash \KK_\actiondom(\a_* \sqrt\vae,\t)\Big)\le \bar c_* \, \torsion_0\ N_\rr^{\rm int}(\DD)\ \mathsf{M}^{-1}\ {\rr^{d-1}}\ \a_*\ \sqrt{\vae}\,,
\eeq
with
\beq{cbstar}
\bar c_*\coloneqq \frac54\, (2\p)^d\, \Big( \frac{d 2^{2d}}{c_0} +2^d d^{\frac{d-1}2}\ \sum_{k\in\integer^d\bks\{0\}}\ \frac1{|k|_1^\t\, |k|_2}\Big)
\,.
\eeq
\ethm

\proof  Let $\a\coloneqq \a_* \sqrt\vae$. Then,   $\r=\mathsf{L} \a/c_0$, so that the first inequality in \equ{smcEAr0v2} is satisfied (with the equal sign). Furthermore, with the above positions, $\epsilon$ in \equ{defpar} is given by 
$$\epsilon=\frac{\mathsf{M}\mathsf{P}}{\a_*^2}\,,
$$
so that also the second inequality in \equ{smcEAr0v2} is satisfied (with the equal sign). Finally, the relation $\r<\rr$ is equivalent to $\vae<\vae_*$, which is satisfied by hypothesis. Hence, all the assumptions of Theorem~\ref{arnoldtheorem} are satisfied
and therefore the measure estimate \equ{MesChPin} holds with $\KK$ as in \equ{KolmogorovSet}.  

\nl
We proceed to estimate the two terms in the right hand side of \equ{MesChPin} separately. Let us first discuss the measure of $\cR_\a$.

\nl 
We claim that {\sl  the map $y\in \rball_\rr(y_0)\mapsto \o_0(y)$ is a diffeomorphism for every $y_0\in \DD$}.
\\
To see this, we shall apply the quantitative Inverse Function Theorem~\ref{IFT} to $f(y)=K_y(y)$. In such a case, we can take $T=K_{yy}(y_0)^{-1}$  
(using Cauchy estimates, see Lemma~\ref{Cau}),
\beqano
\|\uno_d-TK_{yy}(y)\|&\le& \|T\|\|K_{yy}(y_0)-K_{yy}(y)\|\\
&\le & { d^2\mathsf{L}\|\partial_y K_{yy}\|_\rr \, \rr \le d^2\mathsf{L}\,  \frac{\|K_{yy}\|_\ro}{\ro-\rr}\, \rr}\\
&\le & d^2\, \mathsf{L}\mathsf{M}\ \frac{\rr}{\ro-\rr}=\su2\,,
\eeqano
where
$\|\dpr_{y}K_{yy}\|_{\rr}\coloneqq \sup_{\rball_{\rr}(y_0)}\max\{|\dpr^3_{y_iy_jy_k}K|:\, i,j,k=1,\cdots,d\}.
$
Hence, by Theorem~\ref{IFT}, $\o_0$ is invertible on any ball $\rball_\rr(y_0)$ with $y_0\in \DD$, as claimed.

\nl
Now, let $\Do=\{y_1,...,y_{n_0}\}$ with 
$n_0\coloneqq N^{\rm int}_\rr(\DD)$. Then:
\beqano
\meas (\cR_\a) &\le & \sum_{j=1}^{n_0} \meas(\{y\in \rball_\rr(y_j):\ \o_0(y)\notin\Dioph_\a^\t\})\nonumber\\
&\le& 
 \sum_{j=1}^{n_0} \sum_{k\in\integer^d\bks\{0\}} \meas\Big(\Big\{y\in \rball_\rr(y_j):\ | \o_0(y)\cdot e_k|\le\frac{\a}{|k|_1^\t|k|_2}\Big\}\Big)\,.
\eeqano
where $e_k\coloneqq\frac{k}{|k|_2}$. 
Since on $\rball_\rr(y_j)$, $y\to\o_0(y)$ is a diffeomorphism, by the change of variables
$y=\o_0^{-1}(\o)$, we find
\beqano
&& \meas\Big(\Big\{y\in \rball_\rr(y_j):\ | \o_0(y)\cdot e_k|\le\frac{\a}{|k|_1^\t|k|_2}\Big\}\Big)\\
&\le& 
\d_0^{-1} \meas\Big(\Big\{\o\in \o_0\big(\rball_\rr(y_j)\big):\ | \o\cdot e_k|\le\frac{\a}{|k|_1^\t|k|_2}\Big\}\Big)\\
&\le & \d_0^{-1} \Big(\diam \o_0\big(\rball_\rr(y_j)\Big)^{d-1} \frac{2\a}{|k|_1^\t|k|_2}\\
&\le & \d_0^{-1} \big(\mathsf{M} 2\sqrt{d}\rr\big)^{d-1}  \frac{2\a}{|k|_1^\t|k|_2}\,.
\eeqano
Summing up over $j$ and $k$, one gets 
\beq{diophest}
\meas (\cR_\a)\le \Big(2^d d^{\frac{d-1}2}\ \sum_{k\in\integer^d\bks\{0\}}\ \frac1{|k|_1^\t\, |k|_2}\Big) \, n_0\, \d_0^{-1}\, \mathsf{M}^{d-1}\,\rr^{d-1}\, \a\,.
\eeq
Let us turn to the estimate of  $\meas({\mathscr T}_\r)$. Observing that
$$
{\mathscr T}_\r\coloneqq \rball_{\rr+\r} (\Do)\backslash  \rball_{\rr-\r} (\Do)\subset \bigcup_{j=1}^{n_0} 
 \rball_{\rr+\r} (y_j)\backslash  \rball_{\rr-\r} (y_j)\,,
$$
one finds
\beqa{measFr}
\meas({\mathscr T}_\r)&\le& \sum_{j=1}^{n_0} \meas\Big(\rball_{\rr+\r} (y_j)\backslash  \rball_{\rr-\r} (y_j)\Big)\nonumber\\
&=& n_0 2^d \big((\rr+\r)^d-(\rr-\r)^d\big)\nonumber\\
&\le& d 2^{2d}\, n_0\, \r\, \rr^{d-1}=\frac{ d 2^{2d}}{c_0} \, n_0\,  \frac{\torsion}{\mathsf{M}}\, \rr^{d-1} \a
\,.
\eeqa
Observing that $\KK\subset \KK_\actiondom(\a_* \sqrt\vae,\t)$ and that,  by \equ{NormGstrThtv21}, $\b$ in \equ{defsets} satisfies $\b<\frac54 (2\p)^d$,  one sees that \equ{diophest} and \equ{measFr} imply \equ{eas1} with $\bar c_*$ as in \equ{cbstar}. \qed

\subsection{Smooth case}\label{sec:sm}
In order to state the result for  case (B), we need the following definitions.

\begin{itemize}

\item[\tiny $\bullet$] Let $S$ be a compact and connected $C^2$--hypersurface of $\rn$.
The     {\bf minimal focal distance of $S$} is defined as: 
$$
{\minfoc({S})\coloneqq \min\big\{ \inf\{e_c(u,\nu^+(u)):u\in S\}\,,\inf  \{e_c(u,\nu^-(u)):u\in S\}\big\}\,,}
$$ 
where $\nu^\pm(u)$ denotes  the outwards/inwards normal to $S$ at $u$ and 
$$
e_c(u,v)\coloneqq \dst\sup \{t>0:\dist_2(u+tv,{S})=t\}\;,
$$
$\dist_2$ being  euclidean distance.

\item[\tiny $\bullet$] Given any bounded set $D$  in $\rn$, we define the (measure of the) {\bf maximal $(d-1)$--dimensional section of $D$} as 
$$\sec_{d-1}(D)\coloneqq \sup_{\l\in \Lambda^{d-1}}
 \cH^{d-1}(\l\cap D)$$
 where $\Lambda^{d-1}$ denotes the set of all hyperplanes in $\rn$ and $\cH^{d-1}$ the $(d-1)$--dimensional Hausdorff measure.
  
\item[\tiny $\bullet$] 
Given a set $D\subset \rn$ and $\r>0$, we  define {\bf $\r$--inner domains of $D$} (which depend upon the choice of  the metric) as\footnote{Recall that $\rball_\r$ denotes a ball with respect to  the sup norm $|\cdot|=|\cdot|_\io$, while ${\bf B}_\r$ denotes a ball with respect to the Euclidean norm $|\cdot|_2$.} 
\beq{innerdom}
D'_\r\coloneqq \big\{y\in D:\ \rball_\r(y)\subset D\big\}\,,\qquad D''_\r\coloneqq \big\{y\in D:\ {\bf B}_\r(y)\subset D\big\}\,,
\eeq

\end{itemize}

\thm{thmmeas3} 
Let $\actiondom\subset \rn$ be a open and bounded set with $C^2$, compact and connected  boundary. Let $\t>d-1\ge1$, $s>0$. Let 
 $K\in \mathcal{B}_\ro(\actiondom)$ be  uniformly (Kolmogorov) non--degenerate and so that the unperturbed frequency map $y\in
 \actiondom\mapsto \o_0(y)\coloneqq K_y(y)\in \complex^d$ is a global diffeomorphism. Let
 $P\in\mathcal{B}_{\ro,s}(\actiondom)$ and let $\mathsf{M}$, $\mathsf{L}$, $\mathsf{P}$, $\torsion$ as in \equ{defpar} and define
\beq{defrrs}
\rr\coloneqq \min\{\ro\,,\,   \minfoc(\partial \actiondom)\,,\, 1/\kappa \}/\sqrt{d}\,,
\eeq
where $\kappa\coloneqq \sup_{\partial \actiondom}\max_{ 1\le j\le d-1}|\kappa_j|$, $\kappa_j$'s being the principal curvatures of $\partial \actiondom$.
 Let $c_0$ and $c_*$ be as in Theorem~\ref{arnoldtheorem};
 fix $0<s_*<s$, let $\epsilon_*$ be as in \equ{smcEAr0v2}; let $\a_*$, $\vae_*$, $\d_0$, $\torsion_0$ and $\r$ 
 be as in \equ{def.meas1}.
Let 
$\KK_\actiondom(\a_* \sqrt\vae,\t)$ be as in \equ{KAMtori} with $H=K+\vae P$.
Then, if $0< \vae<\vae_*$, 
one has
\beq{eas2}
\meas\Big((\actiondom\times\torus^d)\backslash \KK_\actiondom(\a_* \sqrt\vae,\t)\Big)\le \hat c_* \, \torsion_0\  
\mathsf{M}^{-1}\ \max\big\{\sec_{d-1}(\actiondom)\,,\, \cH^{d-1}(\partial \actiondom)\big\}\ 
\a_*\ \sqrt{\vae}\,,
\eeq
with
\beq{chstar}
\hat c_*\coloneqq \frac52\, (2\p)^d\, \Big( \frac{2^{d}}{\sqrt{d}\ c_0} +  \sum_{k\in\integer^d\bks\{0\}}\ \frac1{|k|_1^\t\, |k|_2}\Big)
\,.
\eeq
\ethm
\proof  The idea is again to apply Theorem~\ref{arnoldtheorem} and Theorem~\ref{thmmeas1}. 
\\
Let ${\Do\coloneqq \actiondom''_{\sqrt{d}\rr}}$.
Since $\sqrt{d}\rr\le \minfoc(\partial\actiondom)$,  by Lemma~\ref{lem:geo1}, 
\beq{yESeQ}
{\rball_\rr(\Do)\subset {\bf B}_{\sqrt{d}\rr}(\actiondom''_{\sqrt{d}\rr})=\actiondom\,,\qquad {\rm and}\qquad  \hactiondom=\rball_{\rr-\r}(\Do)\supseteq{\bf B}_{\rr-\r}(\Do)=\actiondom''_{(\sqrt{d}-1)\rr+\r}}\,. 
\eeq
As in the proof of Theorem~\ref{thmmeas2}, we
let $\a\coloneqq \a_* \sqrt\vae$, so that   $\r=\mathsf{L} \a/c_0$ and  $\epsilon=\mathsf{M}\mathsf{P}/\a_*^2$ (cfr.  \equ{defpar}).
Then,   the  inequalities in \equ{smcEAr0v2}  
hold with the equal sign.
The relation $\r<\rr$ is equivalent to $\vae<\vae_*$, which is satisfied by hypothesis. Hence, all the assumptions of Theorem~\ref{arnoldtheorem} are satisfied
and  the measure estimate \equ{MesChPin} holds with $\KK$ as in \equ{KolmogorovSet}.  

\nl
By hypothesis the frequency map $y\to\o_0(y)$ is a diffeomorphism on $\actiondom$, so we can repeat the estimate 
on the   measure of $\cR_\a$ done in the proof of Theorem~\ref{thmmeas2} without the need of localizing the actions. Letting, as above, $e_k\coloneqq\frac{k}{|k|_2}$, we find
\beqano
\meas (\cR_\a) &=&   \meas(\{y\in \actiondom :\ \o_0(y)\notin\Dioph_\a^\t\})\nonumber\\
&\le& 
\sum_{k\in\integer^d\bks\{0\}} \meas\Big(\Big\{y\in \actiondom:\ | \o_0(y)\cdot e_k|\le\frac{\a}{|k|_1^\t|k|_2}\Big\}\Big)\,.\\
&\le& 
\sum_{k\in\integer^d\bks\{0\}}\d_0^{-1} \meas\Big(\Big\{\o\in \o_0\big(\actiondom)\big):\ | \o\cdot e_k|\le\frac{\a}{|k|_1^\t|k|_2}\Big\}\Big)\\
&\le & \d_0^{-1} \, \mathsf{M}^{d-1}\, \sec_{d-1}(\actiondom) \sum_{k\in\integer^d\bks\{0\}}  \frac{2\a}{|k|_1^\t|k|_2}\\
&\le & \torsion_0\, \mathsf{M}^{-1}\, \sec_{d-1}(\actiondom) \sum_{k\in\integer^d\bks\{0\}}  \frac{2\a}{|k|_1^\t|k|_2}\,.
\eeqano
The estimate on the measure of ${\mathscr T}_\r$ follows from Lemma~\ref{lem:geo2}. Indeed, 
if we denote  ${ \mathfrak T}_{\r}(S)\coloneqq  \{u\in \rn:\ \dist_2(u,S)< \r\}$, we have
 (compare \equ{def.t.n})
$$
{{\mathscr T}_\r=
\rball_{\rr+\rl}(\Do)\setminus\rball_{\rr-\rl}(\Do)\stackrel{\equ{yESeQ}}{\subseteq}{\bf B}_{\sqrt{d}\r}(\actiondom)\bks\actiondom''_{(\sqrt{d}-1)\rr+\r} \subseteq {\mathfrak T}_{\sqrt{d}\rr}(\partial\actiondom)\,,}
$$
Since $\rr\le \min\{\minfoc(\partial \actiondom)/\sqrt{d}\,,\, 1/(\sqrt{d}\kappa) \}$, by \equ{minfoc2}, we get
\beqano
\meas({\mathscr T}_\r)&\le&
\meas ({\mathfrak T}_{\sqrt{d} \rr}(\partial\actiondom))
\\
&\le& \frac2d\, \frac{(1+\sqrt{d}\rr\kappa)^{d}-1}{\kappa}\,\mathcal{H}^{d-1}(\partial\actiondom)\\
&<& \frac{2^{d+1}}{\sqrt{d}}\,\r\, \mathcal{H}^{d-1}(\partial\actiondom)\\
&=& \frac{2^{d+1}}{\sqrt{d}c_0}\,\torsion\, {\mathsf M}^{-1} \, \mathcal{H}^{d-1}(\partial\actiondom)\, \a\,.
\eeqano
Since $\a=\a_*\sqrt{\vae}$, \equ{eas2} follows, with $\hat c_*$ as in \equ{chstar}. \qed

\appendix
\section*{Appendix}
\addcontentsline{toc}{section}{Appendices}
\setcounter{section}{0}
\renewcommand{\thesection}{\Alph{section}}

\appA{Proof of Theorem~\ref{arnoldtheorem}}
In this appendix we provide the details needed to prove  Arnold's Global KAM Theorem (Theorem 1). The main point is the choice of the various parameters and sequences involved in the Newton--like procedure based on the iteration of a `KAM step'  
(in turn, based upon the original scheme by Arnold; compare \cite{ARV63} and its revisions in \cite{koudjinan2019quantitative} and \cite{chierchia2019ArnoldKAM}).  Although the main ideas are well known, some details are needed, especially in order to compute explicitly constants and to keep the optimal relation between $\vae$ and $\a$. Furthermore, the construction of the `integrating map' also require a discussion. All this is done in the present appendix.

\nl
By following \cite[Chap.~6]{koudjinan2019quantitative}, one gets the following:

\subsection*{General step of the KAM scheme}
\lemtwo{lem:1Extv5} {KAM step}
Let $r>0,\; 0<2\s< s\le 1$, $\actiondom_\sharp\subset\rn$ be a non--empty, bounded domain. Consider the Hamiltonian parametrized by $\vae\in \real$
$$
H(y,x;\vae)\coloneqq K(y)+\vae P(y,x)\;,
$$
where $K,P\in \mathcal{\rball}_{r,s}(\actiondom_\sharp)$. 
Assume that\footnote{In the sequel, $K$ and $P$  stand for  generic real analytic Hamiltonians which, later on, will respectively play the roles of $K_j$ and $P_j$,  and $y_0,\,r$, the roles of $y_j,\,r_j$ in the iterative step.} 
\beq{RecHypArnExtv5}
\begin{aligned}
&
\det K_{yy}(y)\not=0
\;,\qquad\qquad\qquad T(y)\coloneqq K_{yy}(y)^{-1}\;,\quad \forall\;y\in\actiondom_\sharp\;,\\
&\|K_{yy}\|_{r,\actiondom_\sharp}\le \mathsf{M}\;,\qquad\qquad\qquad\ \, \|T\|_{\actiondom_\sharp}\le \mathsf{L}\;,\\
& \|P\|_{r,s,\actiondom_\sharp}\le \mathsf{P} \;,\qquad\qquad\qquad\,\,\, K_y(\actiondom_\sharp)\subset \D^\t_\a\;. 
\end{aligned}
\eeq
\noi
Fix  $\vae\neq0$ and assume that
\beq{DefNArnExt1v5}
\l\ge \log\left(\s^{2\n+d}\frac{{\a}^2}{ {\vae}{\mathsf{P}}\mathsf{M}}\right)\ge 1 \;.
\eeq
Let
\beq{DefNArn2v5}
\begin{aligned}
&
\k\coloneqq 4\s^{-1}\l\;, \quad\ \ \, \check{r}\le \frac{r}{32d\mathsf{L} \mathsf{M}}\;, \quad 
\bar{r}\le
\dst\min\left\{\frac{\a}{2d\mathsf{M}\k^\n}\,,\, \check{r} \right\},\\
&\tilde r\coloneqq \frac{\check{r}\s}{16d\mathsf{L} \mathsf{M}},\quad\quad \bar{s}\coloneqq s-\frac{2}{3}\s,\quad\, s'\coloneqq s-\s \,,
\end{aligned}
\eeq
and\footnote{Notice that ${\smalpara}\ge \s^{-d}\ovl{{\smalpara}}\ge \ovl{{\smalpara}}$ since $\s\le 1$. Notice also that $\mathsf{L}\mathsf{M}\ge 1$, so that $\frac{16\mathsf{L} }{r\check{r}}\s^{-(\n+d)}> \frac{16\mathsf{L} }{r^2}\ge \frac{4}{\mathsf{M}r^2}$.
%$$\ovl{{\smalpara}}=\mathsf{C}_0 \left(\frac{r_0}{\a}\right)^2\frac{\mathsf{M}_\infty \mathsf{P}}{r_0 r}\s^{-(2\t+d+1)}\le \mathsf{C}\frac{\mathsf{L}_\infty^2\mathsf{M}_\infty \mathsf{P}}{r_0 r}\s^{-2(\t+d+1)}\le {\smalpara}\,.$$
}
\begin{align*}
{\smalpara}&\coloneqq \mathsf{P}\dst\max\left\{\frac{16\mathsf{L}  }{r\bar{r}}\s^{-(\n+d)}\,,\,
 \frac{\mathsf{C}_4 }{\a\bar{r} }\s^{-2(\n+d)}\right\}\;.
\end{align*}
Assume:   
\beq{cond1ExtExtv5}
 {\vae}{{\smalpara}}\le \frac{\sigma}{3}
\ .
\eeq
Then, there exists 
a  diffeomorphism $G\colon  \cball_{\tilde{r}}(\actiondom_\sharp){\to}G( \cball_{\tilde{r}}(\actiondom_\sharp))$, a symplectic change of coordinates
\beq{phiokExt0Ext}
\phi'=\id+\vae \tilde{\phi}: \cball_{\bar{r}/2,s'}(\actiondom_\sharp')\to \cball_{2r/3, \bar{s}}(\actiondom_\sharp),
\eeq
such that 
\beq{HPhiH'}
\left\{
\begin{aligned}
& H\circ \phi'\eqqcolon H'\eqqcolon K'+\vae^2 P'\ ,\\
& \dpr_{y'}K'\circ G=\dpr_{y}K,\quad \det \dpr_{y'}^2K'\circ G\neq 0\quad\mbox{ on } \actiondom_\sharp\,,
\end{aligned}
\right.
\eeq
with $K'(y')\coloneqq K(y')+\vae\wt K(y')\coloneqq K(y')+\vae \average{P(y',\cdot)}$. Moreover, letting $\left(\dpr^2_{y'} K'(\mathsf{y}')\right)^{-1}\eqqcolon T(\mathsf{y}')+\vae\;\wt T(\mathsf{y}')$, $\mathsf{y}'\in G(\actiondom_\sharp)$, the following estimates hold.

\beq{convEstExt}
\left\{
\begin{aligned}
&\|\dpr_{y'}^2\wt K\|_{r/2,\actiondom_\sharp}\le  
\mathsf{M}{\smalpara}\,,\qquad \|G-\id\|_{\tilde{r},\actiondom_\sharp}\le \s^{\n+d}\bar{r} {\vae}{\smalpara}\;,\qquad \|\wt T\|_{\actiondom_\sharp'}\le \mathsf{L}{\smalpara}\,,\\
& \max\left\{\frac{ \mathsf{C}_{12} }{\mathsf{C}_4}\|\ovl{\mathsf{W}}\nabla\tilde \phi\;\ovl{\mathsf{W}}^{-1}\|_{\bar{r}/2,s',\actiondom_\sharp'},\|\mathsf{W}\,\tilde \phi\|_{\bar{r}/2,s',\actiondom_\sharp'}\right\}\le \s^{d}{{\smalpara}}\,,\qquad \|P'\|_{\bar{r}/2, s',\actiondom_\sharp'} \le  {\smalpara}\mathsf{P}\,, 
\end{aligned}
\right. 
\eeq
where
\begin{equation*}
\actiondom_\sharp'\coloneqq G(\actiondom_\sharp)\;, \quad \left(\dpr^2_{y'} K'(\mathsf{y}')\right)^{-1}\eqqcolon T\circ G^{-1}(\mathsf{y}')+\vae\;\wt T(\mathsf{y}')\,, \ \forall\; \mathsf{y}'\in \actiondom_\sharp'\,,
\end{equation*}
$$
\mathsf{W}\coloneqq \diag(\bar{r}^{-1}\uno_d,\uno_d)\;,\qquad \ovl{\mathsf{W}}\coloneqq  \diag(\s^{-\t}\bar{r}^{-1}\uno_d,\uno_d)\;.
$$
\elem
\subsubsection*{Implementation}
As in \cite{chierchia2019ArnoldKAM}, we shall separate the first step from the others. Let   $H$, $K$, $P$, $\r$, $s$, $s_*$,  $\mathsf{W}$, ${\mathsf{P}}$, $\mathsf{M}$, $\mathsf{L}$, $\torsion$, $\epsilon$ be  as in \S 2. Set 
\begin{align*}
&\s_0\coloneqq (s-s_*)/2\;,\ \epsilon_0\coloneqq \epsilon\;,\ \torsion_0\coloneqq \torsion\;,\ r_0\coloneqq \r\;,\ \mathsf{L}_0\coloneqq \mathsf{L}\;,\ \mathsf{M}_0\coloneqq \mathsf{M}\;,\ \mathsf{P}_0\coloneqq \mathsf{P}\;,\ \mathsf{W}_1\coloneqq \mathsf{W}\;,\\
&\l_0\coloneqq \log\epsilon_0^{-1}\;,\  \l_*\coloneqq \mathsf{C}_7 \s_0^{-(4\n+2d+1)}\torsion_0^2\l_0^{2\n}\;,\   \torsion_*\coloneqq 2^{2\n+2d+1}\mathsf{C}_5^2\torsion_0^2\;,\  \k_0 \coloneqq 4\s_0^{-1}\l_0\;,\\
&K_0\coloneqq K\;,\ P_0\coloneqq P\;,\ H_0\coloneqq H\;,\quad  \actiondom_0\coloneqq \Can\;.
\end{align*}
\subsubsection*{First Step}
\noi
Let
\begin{align*}
&s_1\coloneqq s_0-\s_0\;,\quad \check{r}_{1}\coloneqq \frac{r_0}{64d\torsion_0}\;, \quad \tilde r_{1}\coloneqq \frac{\check{r}_{1}\s_0}{32d\torsion_0}\;,\quad r_1\coloneqq\su2\min\left\{\frac{\a}{2d\sqrt{2}\mathsf{M}_0\k_0^{\n}}\,,\, \check{r}_1 \right\}\;,\\
& \mathsf{M}_1\coloneqq \left(1+\frac{\s_0}{3}\right)\mathsf{M}_0 \;,\quad \mathsf{L}_1\coloneqq \left(1+\frac{\s_0}{3}\right)\mathsf{L}_0\;,\quad \hat{\epsilon}_0\coloneqq \mathsf{C}_{8}\s_0^{-(3\n+2d+1)}\epsilon_0^{1/2},\quad  \mathsf{P}_1\coloneqq \frac{\hat{\epsilon}_0 \mathsf{P}_0}{ {\vae}},\\
&{\smalpara}_0\coloneqq \mathsf{P}_0\dst\max\left\{\frac{8\mathsf{L}_0 }{r_0 r_1}\s_0^{-(\n+d)}\,,\, \frac{\mathsf{C}_4}{2\a r_1}\s_0^{-2(\n+d)}\right\}\;.
\end{align*}
\lem{frstStepv5}
Under the above assumptions and notations, if
\beq{condBisv2v5}
\a\le \frac{\mathsf{C}_4}{16}\frac{r_0}{\mathsf{L}_0}\qquad\mbox{and}\qquad \max\left\{\ex\;\epsilon_0\;,\, \hat{\epsilon}_0\right\}\le 1\;,
\eeq
then, there exist $\actiondom_1\subset \actiondom$, a real--analytic diffeomorphism
$$
G_{1}\colon  \cball_{\tilde{r}_{1}}(\Can){\to}G_{1}( \cball_{\tilde{r}_{1}}(\Can))
$$
and a real--analytic symplectomorphism 
\beq{phijExtv5}
\phi_{1}:\cball_{r_{1},s_{1}}(\actiondom_1)\to \cball_{r_{0},s_{0}}(\actiondom_{0})
\eeq
such that
\begin{align}
&G_1(\Can)=\actiondom_1\;,\label{k1k0puv5}\\
&\dpr_{y_1}K_{1}\circ G_{1}=\dpr_{y}K_0 \;,\label{kjkjpuv5}\\
&H_{1}\coloneqq H_{0}\circ\phi_{1}\eqqcolon K_{1} + \vae^{2} P_{1}\qquad\quad\ \quad\; \mbox{on } \cball_{r_{1},s_{1}}(\actiondom_{1})\label{HjExtExtv5}
\end{align}
and\footnote{\equ{estG1devidv2} follows trivially \equ{estG1idv2} using Cauchy's estimate.}
\begin{align}
&\actiondom_1\subset \actiondom_{r_1}\;,\label{estfin2Bis0001v2}\\
& \|\dpr_{y_1}^2 K_1\|_{r_0/4,\actiondom_1}\le \mathsf{M}_1\;,\qquad \|T_1\|_{\actiondom_1}\le  \mathsf{L}_1\;,\qquad T_1\coloneqq (\dpr_{y_1}^2 K_1)^{-1}\;,\label{estfin2Bis0000v2}\\
&  \|P_1\|_{r_1,s_1,\actiondom_1}\le  \mathsf{P}_1\;,\label{estfin2Bis000v2}\\
&
\|G_1-\id\|_{\tilde r_1,\Can}\le 2\s_0^{\n+d}\;r_1\; {\vae}{\smalpara}_0\;,\label{estG1idv2}\\
&\|\dpr_z G_{1}-\uno_d\|_{\tilde{r}_{1}/2,\Can}\le 2^{5}d\mathsf{C}_4\sqrt{2}\torsion_0\s_0^{\t+d}\k_0^{-\n} {\vae}{\smalpara}_0 \;,\label{estG1devidv2}\\
&
\max\{ \mathsf{C}_{12} \mathsf{C}_4^{-1}\|\ovl{\mathsf{W}}_1\nabla(\phi_1-id)\ovl{\mathsf{W}}_1^{-1}\|_{r_1,s_1,\actiondom_1},\;\|\mathsf{W}_1(\phi_1-\id)\|_{r_1,s_1,\actiondom_1}\}
\le \s_0^d\; {\vae}{{\smalpara}}_0\;.\label{estfin2Bis010v2}
\end{align}
\elem
\subsubsection*{Second step, iteration and convergence}
For a given $j\ge 1$,  define\footnote{Notice that $s_{j}\downarrow s_*$ and $r_{j}\downarrow 0$.}
\begin{align*}
& \dst\s_j\coloneqq \frac{\s_0}{2^j}\,,\quad s_{j+1}\coloneqq s_j-\s_j=s_*+\frac{\s_0}{2^j}\,,\quad  \bar s_{j}\coloneqq s_j-\frac{2\s_j}{3}\,,\quad \k_j \coloneqq 4^j\k_0\,,\\
&\mathsf{M}_{j+1}\coloneqq \mathsf{M}_0\dst\prod_{k=0}^{j}(1+\frac{\s_k}{3})<\mathsf{M}_0\sqrt{2}\,,  \quad \mathsf{L}_{j+1}\coloneqq \mathsf{L}_0\dst\prod_{k=0}^{j}(1+\frac{\s_k}{3})<\mathsf{L}_0\sqrt{2}\,,\\
& \epsilon_j\coloneqq \frac{\mathsf{M}_0 {\vae}^{2^j} \mathsf{P}_j}{ \a^2}\,,\quad \check{r}_{j+1}\coloneqq \frac{r_j}{64d\torsion_0}\,,\quad \tilde r_{j+1}\coloneqq \frac{\check{r}_{j+1}\s_j}{32d\torsion_0}\,,\quad r_{j+1}\coloneqq \su2\min\left\{\frac{\a}{2d\sqrt{2}\mathsf{M}_0\k_{j}^{\n}}\,,\, \frac{r_{j}}{64d\torsion_0} \right\},\\
& \mathsf{P}_{j+1}\coloneqq 	\l_* \torsion_*^{j-1} \frac{\mathsf{M}_0 \mathsf{P}_{j}^2}{\a^2}\,,\quad \hat{\epsilon}_j     \coloneqq \l_* \torsion_*^j\epsilon_j\,,\quad \mathsf{W}_{j+1}\coloneqq \diag\left((2r_{j+1})^{-1}\uno_d\,,\uno_d\right)\,,\\ 
&\ovl{\mathsf{W}}_{j+1}\coloneqq \diag\left(\s_j^{-\t}(2r_{j+1})^{-1}\uno_d\,,\uno_d\right)\,,\quad {\smalpara}_j\coloneqq \mathsf{P}_j\dst\max\left\{ \frac{8\mathsf{L}_0\sqrt{2} }{r_jr_{j+1}}\s_j^{-(\n+d)}\,,\,
\frac{\mathsf{C}_4 }{2\a r_{j+1} }\s_j^{-2(\n+d)}\right\}\,.
\end{align*}
Observe that, for any $j\ge1$,
\begin{align*}
\hat{\epsilon}_{j+1}&= \l_* \torsion_*^{j+1}\epsilon_{j+1}=\l_* \torsion_*^{j+1}\frac{\mathsf{M}_0 {\vae}^{2^{j+1}} \mathsf{P}_{j+1}}{ \a^2}=\l_* \torsion_*^{j+1}\frac{\mathsf{M}_0  {\vae}^{2^{j+1}}}{ \a^2}\;\l_* \torsion_*^{j-1} \frac{\mathsf{M}_0 \mathsf{P}_j^2}{\a^2} \left(\l_* \;\torsion_*^{j}\;\epsilon_j\right)^2=\hat{\epsilon}_j^2
\end{align*}
\ie
$$
\hat{\epsilon}_j=\hat{\epsilon}_1^{2^{j-1}} \;.
$$
\lem{lem:2Extv5} 
Assume $\equ{HjExtExtv5}\div\equ{estfin2Bis000v2}$ with some $\vae\neq 0$ and
\beq{condBisv2Prtv5} 
\max\left\{\ex\;\epsilon_0\;,\,2^{11}d^{2}\torsion_0\s_0^{\n+d}\hat{\epsilon}/3\;,\, 2\mathsf{C}_{6}\torsion_0\hat{\epsilon}_{1}\right\}\le 1\;.
\eeq
Then, one can construct a sequence of real--analytic diffeomorphisms 
$$
G_{j}\colon  \cball_{\tilde{r}_{j}}(\actiondom_{j-1}){\to}G_{j}( \cball_{\tilde{r}_{j}}(\actiondom_{j-1}))\;,\qquad j\ge 2
$$
and of real--analytic symplectic transformations 
\beq{phijBisv2v5}
\phi_{j}:\cball_{r_{j},s_{j}}(\actiondom_{j})\to \cball_{r_{j-1},s_{j-1}}(\actiondom_{j-1})\;,
\eeq
such that
\begin{align*}
&G_j(\actiondom_{j-1})=\actiondom_j\subset \actiondom_{r_j} \;,\\
&\dpr_{y}K_{j+1}\circ G_{j+1}=\dpr_{y}K_j \;,\\
&H_{j}\coloneqq H_{j-1}\circ\phi_{j}\eqqcolon K_{j} + \vae^{2^{j}} P_{j}\qquad\quad\  \mbox{on}\quad \cball_{r_{j},s_{j}}(\actiondom_{j})\;,
\end{align*}
converge uniformly. More precisely, we have the following:
\begin{itemize}
\item[$(i)$] the sequence $G^{j}\coloneqq G_{j}\circ G_{j-1}\circ\cdots\circ G_2\circ G_1$ converges uniformly on $\Can$ to a lipeomorphism $Y^*\colon \Can\to \Canstar\coloneqq Y^*(\Can)\subset\actiondom$ and $Y^*\in C^\infty_W(\Can)$\;.  
\item[$(ii)$] $\vae^{2^j}\dpr_y^{\b} P_j$ converges uniformly on 
$\Canstar\times\dst\torus^d_{s_*}$ to $0$, for any $\b\in \natural_0^d$\;;
\item[$(iii)$] $\phi^j\coloneqq \phi_2\circ \cdots\circ \phi_j$ converges uniformly on 
$\Canstar\times\tn$ to a symplectic transformation 
$$
\phi^*\colon \Canstar\times\tn\overset{into}{\longrightarrow} \rball_{r_1}(\actiondom_1)\times\tn,
$$
 with $\phi^*\in C^\infty_W(\Canstar\times\tn)$ and $\phi^*(y,\cdot)\colon \torus^d_{s_*}\ni x\mapsto \phi^*(y,x)$ holomorphic, for any $y\in\Canstar$\;;
\item[$(iv)$]  $K_j$ converges uniformly on 
$\Canstar$ to a function $K_*\in C^\infty_W(\Canstar)$, with
\begin{align*}
&\dpr_{y_*}K_*\circ Y^*=\dpr_{y}K_0 \quad \qquad\qquad\quad\quad\mbox{on} \quad \Can\;,\\
 &\dpr_{y_*}^{\b}(H_1\circ\phi^*)(y_*,x)=\dpr_{y_*}^{\b} K_*(y_*)\;,\quad \forall(y_*,x)\in\Canstar\times\tn\;, \forall\;\b\in \natural_0^d\;.
\end{align*}
\end{itemize}
Finally, the following estimates hold for any $i\ge 2$:\footnote{Observe that \equ{estGidevidv2} follows \equ{estGiidv2} using Cauchy's estimate.}
\begin{align}
&\|G_i-\id\|_{\tilde r_{i},\actiondom_{i-1}}\le 2 r_{i}\;\s_{i-1}^{\n+d}\; {\vae}^{2^{i-1}}{\smalpara}_{i-1}\;,\label{estGiidv2}\\
&\|\dpr_z G_{i}-\uno_d\|_{\tilde{r}_{i}/2,\actiondom_{i-1}}\le 2^{5}d\torsion_0\;\s_{i-1}^{\t+d}\; {\vae}^{2^{i-1}}{\smalpara}_{i-1} \;,\label{estGidevidv2}\\
&\|P_i\|_{r_i,s_i,\actiondom_i}\le  \mathsf{P}_i\ ,\label{estfin2Ext01v5}\\
&\|\mathsf{W}_{2}(\phi^{i+2}-\phi^{i+1})\|_{r_{i+2},s_{i+2},\mathscr{D}_{i+2}}\le \mathsf{a}_2\;  \left(\mathsf{C}_6\torsion_0^{\su4} \hat{\epsilon}_{1}\right)^{2^{i}}\;,\label{phini}\\
&
|\mathsf{W}_2(\phi^*-\id)|
%\}
\le \frac{2\s_0^{d+1}\;\hat{\epsilon}_1}{3\; \torsion_*}
\qquad\qquad\ \mbox{on}\quad \Canstar\times\torus^d_{s_*}\label{estfin2Ext03v5}\ ,
\end{align}
where
$$
\mathsf{a}_2\coloneqq \mathbf{a}_1\;\s_2^d\;\|\mathsf{W}_{2}\phi_2\|_{r_{2}, s_2,\mathscr{D}_{2}}\;.
$$ 
\elem 
\Giu
We can now complete the proof of Theorem~\ref{arnoldtheorem}. First of all, observe that
\beq{log1su2}
(\log t)^{a}\le\left(\frac{2a}{\ex}\right)^a \sqrt{t},\qquad \forall\; t\ge \ex,\quad\forall\; a> \su2.
\eeq
and from the proof, we have
\begin{align}
 {\vae} {\smalpara}_0 (3 \sigma_0^{-1})&\overset{\equ{condBisv2v5}}{\le}6d\mathsf{C}_4\sqrt{2}\s_0^{-2(\n+d)-1}\frac{ \mathsf{K}_0 {\vae} \mathsf{P}_0}{\a^2}\k_0^\n\label{IneqL0S1}\\                                   	&\leby{log1su2} \hat{\epsilon}_0\leby{condBisv2v5} 1,\label{L0ifsg3}
\end{align}
and, for $j\ge 1$,
\beq{eqfRacpj}
 {\vae}^{2^j} {\smalpara}_j (3 \sigma_j^{-1})\le {\hat{\epsilon}_1^{2^{j-1}}}/{ \torsion_*}.
 \eeq
 Let $\phi_*\coloneqq \phi_1\circ \phi^*$. Thus, uniformly on $\mathscr{D}_*\times \torus^d_{s_*}$,\footnote{Observe that $\l_0^{2\n}\epsilon_0\leby{condBisv2Prtv5} (4\n)^{2\n}\sqrt{\epsilon_0}\leby{condBisv2Prtv5} (4\n)^{2\n}(2^{11}d^2\mathsf{C}_5)^{-1/2}\torsion_0^{-1}\s_0^{(\n+d+1)/2}$.} 
\begin{align}
|\mathsf{W}_1(\phi_*-\id)|&\le |\mathsf{W}_1(\phi_1\circ \phi^*-\phi^*)|+|\mathsf{W}_1(\phi^*-\id)|\nonumber\\
&\le \|\mathsf{W}_1(\phi_1-\id)\|_{r_1,s_1,\mathscr{D}_1}+\|\mathsf{W}_1\mathsf{W}_2^{-1}\|\;|\mathsf{W}_2(\phi^*-\id)|\nonumber\\
&\le \s_{0}^d\; {\vae}{{\smalpara}}_0+\frac{2\s_0^{d+1}\;\hat{\epsilon}_1}{3\; \torsion_*}\nonumber\\
&\stackrel{\equ{log1su2}+\equ{IneqL0S1}}{\le}6d\mathsf{C}_4\sqrt{2}\s_0^{-(2\n+2d+1)} \epsilon_0\k_0^\n  +\left(\frac{\n}{2\ex}\right)^\n\mathsf{C}_7\s_0^{-(6\n+4d+2)}\torsion_0^2\k_0^{\n} \epsilon_0\nonumber\\
&\le  \mathsf{C}_9 \torsion_0^2\k_0^{\n} \epsilon_0\;,\nonumber
\end{align}
\ie \equ{estArnTrExtv20}. 
 Moreover, setting $G_{0}\coloneqq \id$, we have 
for any $i\ge 3$,
\begin{align*}
\|G^{i}-\id\|_{\Can}&\le \dst\sum_{j=0}^{i-1}\|G^{j+1}-G^{j}\|_{\Can}=\dst\sum_{j=0}^{i-1}\|G_{j}-\id\|_{\actiondom_{j-1}}\overset{\equ{estGiidv2}+\equ{estG1idv2}}{\le}2\dst\sum_{j=0}^{i-1}r_{j+1}\s_j^{\n+d} {\vae}^{2^j}{\smalpara}_j\\
		&\le 2\;r_1\s_0^\n\dst\sum_{j=0}^\infty \s_j^{d} {\vae}^{2^j}{\smalpara}_j\stackrel{\equ{L0ifsg3} +\equ{eqfRacpj}}{\le} 2r_1\s_0^\n\cdot  \mathsf{C}_9 \torsion_0^2\k_0^{\n}\epsilon_0\;,
\end{align*}
and then passing to the limit, we get
$$
\|Y^{*}-\id\|_{\Can}\le  
2^{2\t+1/2}d^{-1} \mathsf{C}_9 \s_0^\n\torsion_0^2\frac{ {\vae} \mathsf{P}_0}{\a}
%=r_*
\;,
$$
 \ie  \equ{NormGstrThtv20}. 
 Now, observing that, for any $j\ge 1$, $\nabla \phi^{j+1}=\nabla \phi^{j}\nabla \phi_{j+1}$,  $\|\ovl{\mathsf{W}}_j\ovl{\mathsf{W}}_{j+1}^{-1}\|=1$ and $\|\ovl{\mathsf{W}}_{j+1}\ovl{\mathsf{W}}_{j}^{-1}\|\le \mathsf{C}_5\;\torsion_0 $, we obtain
\begin{align*}
\|\ovl{\mathsf{W}}_{1}(\nabla\phi^{j+1}-\uno_{2d})\ovl{\mathsf{W}}_{j+1}^{-1} \|_*&\le \left( \|\ovl{\mathsf{W}}_{1}(\nabla\phi^{j}-\uno_{2d})\ovl{\mathsf{W}}_{j}^{-1} \|_*+1\right)\left(\|\ovl{\mathsf{W}}_{j+1}(\nabla\phi_{j+1}-\uno_{2d})\ovl{\mathsf{W}}_{j+1}^{-1} \|_*+1\right)\\
	&\qquad -1\\
	&\le \left( \|\ovl{\mathsf{W}}_{1}(\nabla\phi^{j}-\uno_{2d})\ovl{\mathsf{W}}_{j}^{-1} \|_*+1\right)\left(\frac{\mathsf{C}_4}{ \mathsf{C}_{12} }\s_{j}^d {\vae}^{2^{j}}{{\smalpara}}_{j}+1\right)-1,
\end{align*} 
which iterated yields\footnote{Use: $\ex^t-1\le t\ex^t$, for any $t\ge 0$.}
\begin{align*}
\|\ovl{\mathsf{W}}_{1}(\nabla\phi^{j+1}-\uno_{2d})\ovl{\mathsf{W}}_{j+1}^{-1} \|_*&\le \prod_{j=1}^\infty \left(\frac{\mathsf{C}_4}{ \mathsf{C}_{12} }\s_{j-1}^d {\vae}^{2^{j-1}}{{\smalpara}}_{j-1}+1\right)-1\\
    &\le \exp\left(\sum_{j=1}^\infty \frac{\mathsf{C}_4}{ \mathsf{C}_{12} }\s_{j-1}^d {\vae}^{2^{j-1}}{{\smalpara}}_{j-1}\right)-1\\
	&\stackrel{\equ{L0ifsg3} +\equ{eqfRacpj}}{\le} \exp( \mathsf{C}_{4} \mathsf{C}_9  \mathsf{C}_{12} ^{-1}\torsion_0^2\k_0^{\n}\epsilon_0)-1\\
	&\le \exp((4d)^{-1})\mathsf{C}_{4} \mathsf{C}_9  \mathsf{C}_{12} ^{-1}\torsion_0^2\k_0^{\n}\epsilon_0\stackrel{\equ{smcEAr0v2}+\equ{log1su2}}{\le} \su{4(18d^3+70)\torsion}\;,
\end{align*} 
 and letting $j\to\infty$, we  obtain
 \beq{dprXu*s}
 \|\dpr_x u_*\|_*\le \exp((4d)^{-1})\mathsf{C}_{4} \mathsf{C}_9  \mathsf{C}_{12} ^{-1}\torsion_0^2\k_0^{\n}\epsilon_0\le \su{4(18d^3+70)\torsion}\;.
 \eeq
 \ie \equ{estArnTrExtv21}.

\Giu
\noi
Next, we show that $\Lip_{\mathscr D^*}(Y^*-\id)<1$, which will imply that\footnote{See Proposition II.2 in \cite{zehnder2010lectures}.} $Y^*\colon \mathscr D^*\overset{onto}{\longrightarrow}\mathscr D_*$ is a lipeomorphism. For, observe first that, for any $j\ge 1$, $0<r<\tilde{r}_{j}/2$, $\mathsf{y}_{j-1}\in\mathscr{D}_{j-1}$ and any $y\in \cball_{r}(\mathsf{y}_{j-1})$, we have
$$
|G_j(y)-G_j(\mathsf{y}_{j-1})|\le |(G_j(y)-y)-(G_j(\mathsf{y}_{j-1})-\mathsf{y}_{j-1})|+|y-\mathsf{y}_{j-1}|\overset{\equ{estGidevidv2}+\equ{condBisv2Prtv5}}{\le}\su2|y-\mathsf{y}_{j-1}|+r<2r,
$$
so that
\beq{InclGjse}
G_j(\cball_{r}(\mathscr{D}_{j-1}))\subseteq \cball_{2r}(G_j(\mathscr{D}_{j-1}))=\cball_{2r}(\mathscr{D}_{j}).
\eeq
Thus, as the sequence $\tilde{r}_j$ is strictly decreasing, for any $j\ge k\ge 1$, $G^k$ is well--defined on  $\cball_{2^{-j-1}\tilde{r}_{j+1}}(\mathscr{D}_{0})$ and we have
\beq{InclGjGUnos}
\begin{aligned}
&G^k(\cball_{2^{-j-1}\tilde{r}_{j+1}}(\mathscr{D}_{0}))\overset{\equ{InclGjse}}{\subseteq}G_j\circ \cdots\circ  G_2 (\cball_{2^{-j}\tilde{r}_{j+1}}(\mathscr{D}_{1}))\overset{\equ{InclGjse}}{\subseteq}\cdots\overset{\equ{InclGjse}}{\subseteq} \cball_{2^{k-j-1}\tilde{r}_{j+1}}(\mathscr{D}_{k})\subseteq\\
	&\subseteq \cball_{2^{-1}\tilde{r}_{k+1}}(\mathscr{D}_{k}).
\end{aligned}
\eeq
Therefore, for any $j\ge 2$,
\begin{align*}
\| G^j-\id\|_{L,\cball_{2^{-j-1}\tilde{r}_{j+1}}(\mathscr D^*)}+1&= \| (G_j-\id)\circ G^{j-1}+(G^{j-1}-\id)\|_{L,\cball_{2^{-j-1}\tilde{r}_{j+1}}(\mathscr D_0)}+1\\
  &\le (\| G_j-\id\|_{L,G^{j-1}(\cball_{\frac{\tilde{r}_{j+1}}{2^{j+1}}}(\mathscr D_{0}))}+1)(\| G^{j-1}-\id\|_{L,\cball_{\frac{\tilde{r}_{j+1}}{2^{j+1}}}(\mathscr D_{0})}+1)\\
  &\leby{InclGjGUnos} (\| G_{j}-\id\|_{L,\cball_{\frac{\tilde{r}_{j}}{2}}(\mathscr D_{j-1})}+1)(\| G^{j-1}-\id\|_{L,\cball_{\frac{\tilde{r}_{j}}{2^{j}}}(\mathscr D_{0})}+1)\\
  &= (\|\dpr_z G_j-\uno_d\|_{\tilde{r}_j/2,\mathscr D_{j-1}}+1)(\| G^{j-1}-\id\|_{L,\cball_{\frac{\tilde{r}_{j}}{2^{j}}}(\mathscr D_{0})}+1)\\
  &\overset{\equ{estGidevidv2}+\equ{estG1devidv2}}{\le} (2^{5}d\torsion_0\;\s_{j-1}^{\t+d}\; {\vae}^{2^{j-1}}\mathsf L_{j-1}+1)(\| G^{j-1}-\id\|_{L,\cball_{\frac{\tilde{r}_{j}}{2^{j}}}(\mathscr D^*)}+1)
\end{align*}
which iterated leads to\footnote{Use, again, $\ex^t-1\le t\ex^t\;,\ \forall\; t\ge 0$, and $2^{5}d\mathsf{C}_4\sqrt{2}\torsion_0\s_0^{\t+d}\k_0^{-\n} {\vae}\mathsf L_0\leby{IneqL0S1} 2^{7}d^2\mathsf{C}_4^2\torsion_0\s_0^{-(\n+d+1)}\epsilon_0\ltby{condBisv2Prtv5} (32d)^{-1}$.}
\begin{align}
\| G^j-\uno_d\|_{L,\mathscr D^*} &\le -1+(1+(32d)^{-1})\dst\prod_{i=2}^\infty (2^{5}d\torsion_0\s_{j-1}^{\t+d} {\vae}^{2^{i-1}}\mathsf{L}_{i-1}+1)\nonumber\\
     &\le -1+\exp((32d)^{-1} +2^{5}d\torsion_0 \sum_{i=1}^\infty \s_{i}^{\t+d} {\vae}^{2^{i}}\mathsf{L}_{i})\nonumber\\
     &\le -1+\exp((32d)^{-1}+2^{5}d\torsion_0\s_{1}^{\t}\sum_{i=1}^\infty \s_{i}^{d} {\vae}^{2^{i}}\mathsf{L}_{i})\nonumber\\
     &\le -1+\exp\left((32d)^{-1}+2^{5}d\torsion_0\s_{1}^{\t}\frac{2\s_0^{d+1}\;\hat{\epsilon}_1}{3\; \torsion_*} \right)\nonumber\\
     &\ltby{condBisv2Prtv5}\ -1+\exp\left((32d)^{-1}+(32d)^{-1}\right)\le \ex^{1/(16d)}/(16d)<\su{4d}\;.\label{LipG^jInf1}
\end{align}
Hence, letting $j\to\infty$, we get that $Y^*$ is Lipschitz continuous, with
$
\Lip_{\mathscr D^*}(Y^*-\id)
$
satisfying \equ{NormGstrThtv21} as
$$
2^{5}d\mathsf{C}_4\sqrt{2}\torsion_0\s_0^{\t+d}\k_0^{-\n} {\vae}\mathsf L_0+\sum_{j\ge 2} 2^{5}d\torsion_0\;\s_{j-1}^{\t+d}\; {\vae}^{2^{j-1}}\mathsf L_{j-1}\leby{condBisv2Prtv5} {c_2} \;\torsion^3\;\;(s-s_*)^{-1}\;\frac{\mathsf{M} {\vae}\mathsf{P}}{\a^2}\;\left(\log\frac{\a^2}{\mathsf{M} {\vae}\mathsf{P}}\right)^\n.
$$ 
Next, we show that $\phi_*\in C^\infty_W(\mathscr D_*\times\tn)$. For any $n,j\ge1$, we have
\begin{align*}
\|G^{n+j}-G^j\|_{\mathscr D^*}&\le \dst\sum_{k=j}^{n+j-1}\|G^{k+1}-G^k\|_{\mathscr D^*}\\
		&\leby{estGiidv2}\dst\sum_{k=j}^{n+j-1}2r_{k+1}\s_k^{\n+d} {\vae}^{2^k}\mathsf L_k\\
		&\le 2r_{j+1}\s_j^{\n}\dst\sum_{k\ge 1}\s_k^{d} {\vae}^{2^k}\mathsf L_k\\
		&\le 2r_{j+1}\s_j^{\n}\;\frac{2\s_0^{d+1}\;\hat{\epsilon}_1}{3\; \torsion_*}\\
		&\ltby{condBisv2Prtv5} \s_j^{\n}\;\tilde{r}_{j+1}\;.
\end{align*}
Now, letting $n\to\infty$, we get
\beq{G^*G^jDist}
\|Y^{*}-G^j\|_{\mathscr D^*}<\s_j^{\n}\;\tilde{r}_{j+1}<\frac{\tilde{r}_{j+1}}{4}\;.
\eeq
Hence\footnote{Recall that, by definition, $G^j(\mathscr D^*)=\mathscr D_j$ and $Y^*(\mathscr D^*)=\mathscr D_*$. }, for any $j\ge 1$,
\beq{HaDstr}
\cball_{\frac{\tilde{r}_{j+1}}{4}}\big(G^j(\mathscr D^*))\stackrel{\equ{G^*G^jDist}}{\subset} \cball_{\frac{\tilde{r}_{j+1}}{2}}(\mathscr D_*)\stackrel{\equ{G^*G^jDist}}{\subset} \cball_{\tilde{r}_{j+1}}(\mathscr D_j)\subset \cball_{r_j}(\mathscr D_j)\;.
\eeq
 Therefore, for any $n\ge 1$, we have
\beqano
\dst\sum_{j\ge 3}\|\mathsf{W}_2(\phi^{j}-\phi^{j-1})\|_{\tilde{r}_{j+1}/2,s_j,\mathscr{D}_*}\left(\frac{\tilde{r}_{j+1}}{2}\right)^{-n}&\leby{HaDstr}& (2^{12}d^2\torsion_0^2)^n\dst\sum_{j\ge 3}\|\mathsf{W}_2(\phi^{j}-\phi^{j-1})\|_{r_j,s_j,\mathscr{D}_j} (r_{j}\s_j)^{-n}\\
	    &\overset{\equ{phini}}{\le}& (2^{12}d^2\torsion_0^2\s_1 r_1)^n\mathsf{a}_2 \dst\sum_{j\ge 3} \left(\mathsf{C}_{6}\torsion_0^{\frac{1}{4}} \hat{\epsilon}_{1}\right)^{2^{j-2}} \left(2\mathsf{a}_1\right)^{n(j-1)}\\
&<& +\infty\;,
\eeqano
since, for $j$ sufficiently large,
$$
 \left(\mathsf{C}_{6}\torsion_0^{\frac{1}{4}} \hat{\epsilon}_{1}\right)^{2^{j-1}}\left(2\mathsf{a}_1\right)^{nj}<\left(\sqrt{2}\mathsf{C}_{6}\torsion_0^{\frac{1}{4}} \hat{\epsilon}_{1}\right)^{2^{j-1}}\leby{condBisv2Prtv5} (1/\sqrt{2})^{2^{j-1}}.
$$
Thus, letting $\Phi_j\coloneqq \phi_1\circ \phi^j$ and using the Mean Value theorem, we have 
\begin{align*}
\dst\sum_{j\ge 2}\|\mathsf{W}_2(\Phi_{j}-\Phi_{j-1})\|_{\tilde{r}_{j+1}/2,s_j,\mathscr{D}_*}\left(\frac{\tilde{r}_{j+1}}{2}\right)^{-n}&\le \|\mathsf{W}_2\nabla\phi_1\mathsf{W}_2^{-1}\|_{{r}_{1},s_1,\mathscr{D}_1}\times\\
  &\times\sum_{j\ge 3}\|\mathsf{W}_2(\phi^{j}-\phi^{j-1})\|_{\tilde{r}_{j+1}/2,s_j,\mathscr{D}_*}\left(\frac{\tilde{r}_{j+1}}{2}\right)^{-n}\\
&< \infty\;.
\end{align*}
Consequently, writing
$$
\Phi_{j}=(\Phi_{j}-\Phi_{j-1})+\cdots+(\Phi_{3}-\Phi_{2})\;,\qquad j\ge 2\;,
$$
and invoking Lemma~\ref{Whit1} (see Appendix \ref{appD}), we conclude that 
$\phi_*=\lim \Phi_j\in C^\infty_W(\mathscr D_*\times\tn)$. 

\Giu
Now, we prove $Y^*\in C^\infty_W(\mathscr D^*)$ analogously. For any $j\ge 2$ and $n\ge 1$, we have
$$
G^j=(G^{j}-G^{j-1})+\cdots+(G^{2}-G^{1})\;,
$$
and, thanks to \equ{InclGjse}, 
 $G^{j+1}-G^j$ is well--defined on  $\cball_{2^{-j-2}\tilde{r}_{j+2}}(\mathscr{D}^*)$, for any $j\ge 1$, so that
\begin{align*}
\dst\sum_{j\ge 1}\|G^{j+1}-G^{j}\|_{\frac{\tilde{r}_{j+2}}{2^{j+2}},\mathscr D^*}\left(\frac{\tilde{r}_{j+2}}{2^{j+2}}\right)^{-n}&=\dst\sum_{j\ge 1}({2^{j+2}}{\tilde{r}_{j+2}}^{-1})^{n}\|(G_{j+1}-\id)\circ G^{j}\|_{\frac{\tilde{r}_{j+2}}{2^{j+2}},\mathscr D^*}\\
    &\leby{InclGjGUnos}\dst\sum_{j\ge 2}({2^{j+2}}{\tilde{r}_{j+2}}^{-1})^{n}\|G_{j+1}-\id\|_{\tilde{r}_{j+1},\mathscr D_{j}}\\
	&\le 2\dst\sum_{j\ge 1}({2^{j+2}}{\tilde{r}_{j+2}}^{-1})^{n}r_{j+1}\s_{j}^{\n+d} {\vae}^{2^j}\mathsf{L}_j\\
	&\ltby{condBisv2Prtv5}\infty\;,
\end{align*}
which proves that $Y^*\in C^\infty_W(\mathscr D^*)$.
 
\Giu
Finally, we prove Kolmogorov's non--degeneracy\footnote{See Appendix~\ref{app9}.} of the Kolmogorov tori $\phi_*(\mathscr D_*\times\torus^d)$.
\\
Fix $y_*\in \mathscr D_*$. Let $y_0\coloneqq (Y^*)^{-1}(y_*)$ and
$$
\hat \epsilon \coloneqq \su{4(18d^3+70)\torsion}\;.
$$
Since $\|\dpr_x u_*\|_{*} \stackrel{\equ{estArnTrExtv21}}{\le}\hat \epsilon<1/2$, then the map $x\longmapsto x+u_*(y_*,x)$ is a diffeomorphism of $\tn$.
Letting
$$(\dpr_x(\id+u_*)(y_*,x))^{-1} \eqqcolon \uno_d+A(y_*,x)\ ,
$$
we have 
\beq{equsttiInv}
\|A\|_{*}\le 2\|\dpr_x u_*\|_{*}\leby{estArnTrExtv21} 2\hat \epsilon< 1\;;\quad
\|v_*\|_{*}\leby{estArnTrExtv21} \frac{ \mathsf{C}_9 \sqrt{2}}{4d}\torsion^2\frac{ {\vae}\mathsf{P}}{\a}\leby{smcEAr0v2} \frac{\mathsf{C}_{4} \mathsf{C}_9 \sqrt{2}}{2^5d\mathsf{C}_{*}}\;\r< \frac{\r}{8}.
\eeq
Moreover, write $K_{yy}(y_*)=K_{yy}(y_0)(\uno_d+K_{yy}(y_0)^{-1}(K_{yy}(y_*)-K_{yy}(y_0)))$ and observe
$$
\dist(y_0,\dpr\mathscr D)\ge\r\quad \mbox{and}\quad |y_*-y_0|\stackrel{\equ{NormGstrThtv20}+\equ{smcEAr0v2}}{\le} \frac{2^{\t-5}\mathsf{C}_{4} \mathsf{C}_9 \sqrt{2}}{d\mathsf{C}_{*}}\;\r<\frac{\r}{64d},
$$
so that
\beq{y0y*KolmNoN}
\dist(y_*,\dpr\mathscr D)\ge \frac{\r}{2}.
\eeq
Thus, by the Mean Value Theorem, we have
\begin{align*}
\|K_{yy}(y_0)^{-1}(K_{yy}(y_*)-K_{yy}(y_0))&\|\leby{y0y*KolmNoN} \mathsf{T}\frac{d^2\mathsf{K}}{{\r/2}}|y_*-y_0|\\
	&{\leby{NormGstrThtv20} 2^{\t+11/2}d^2 \mathsf{C}_9 \;\torsion^3\;\frac{\mathsf{K} {\vae} \mathsf{P}_0}{\a \r}\leby{smcEAr0v2} \frac{2^{\t+15/2}d^2\mathsf{C}_{4} \mathsf{C}_9 }{\mathsf{C}_{*}}}\le \su2.
\end{align*}
Hence, $K_{yy}(y_*)$ is invertible and $\|K_{yy}(y_*)^{-1}\|\le 2\|K_{yy}(y_0)^{-1}\|\le 2\mathsf{T}.$

\noi
In \cite{salamon2004kolmogorov} it is proven that the map
$$
\phi^{y_*}(y,x)\coloneqq (y_*+v_*(y_*,x)+y+ A^{T}y,x+u_*(y_*,x)). 
$$
is symplectic. Then, 
$$
H\circ \phi^{y_*} (y,x)= E^{y_*}+\o^{y_*} \cdot y + Q^{y_*}(y,x) 
$$
with: 
\begin{align*}
&E^{y_*}= K(y_*),\quad \o^{y_*}\coloneqq K_{y}(y_0),\quad \langle Q_{yy}^{y_*}(0,\cdot)\rangle=K_{yy}(y_*)+\average{\mathcal{M}}\;,\\
& \mathcal{M}\coloneqq \dpr^2_y\bigg(K(y_*+v_*+y+A^{T}y)-\su2 y^TK_{yy}(y_*)y\bigg)\Big|_{y=0}+\dpr^2_y(\vae P\circ\phi)\Big|_{y=0}\;,\\
& \|K_{yy}(y_*)^{-1}\mathcal{M}\|_*\le 2\mathsf{T}\mathcal{M}\stackrel{\equ{equsttiInv}}{\le} 2(18d^3+70)\hat \epsilon\torsion= 1/2,
\end{align*}
which show that $\langle Q_{yy}^{y_*}(0,\cdot)\rangle$ is invertible.
\qed

\rem{AssumpExtArnolv2}
Here we list all the constants, which appear in the above proof  and give the explicit expression for the constants $c_k$'s appearing in the statement of Theorem~\ref{arnoldtheorem}. \\
Recall that 
$\t> d-1\ge 1$ and notice that all the $\mathsf{C}_i$'s are  greater than $1$ and depend only upon $d$ and $\t$.
\begin{align*}
&\n	\coloneqq \t+1,\quad \mathsf{C}_0 \coloneqq 4\left(\frac{3}{2}\right)^{2\n+d}\dst\int_{\rn} \left( |y|_1^{\n}+d|y|_1^{2\n}\right)\ex^{-|y|_1}dy,\quad \mathsf{C}_1 \coloneqq 2\left(\frac{3}{2}\right)^{\n+d}\dst\int_{\rn} |y|_1^{\n}\ex^{-|y|_1}dy,\\
&\mathsf{C}_2 \coloneqq 2^{3d}d,\quad \mathsf{C}_3 \coloneqq	d^2\mathsf{C}_1^2+6d\mathsf{C}_1 +\mathsf{C}_2,\quad \mathsf{C}_4 \coloneqq \max\left\{(1+d^2)\mathsf{C}_0,\mathsf{C}_3\right\},\quad \mathsf{C}_5\coloneqq \max\left\{2^{2\n}\,,\,2^7d\right\},\\
&\mathsf{C}_6 \coloneqq \left(2^{-d}\mathsf{C}_5\right)^{\su4},\quad \mathsf{C}_7   \coloneqq 3\cdot 2^{4\n+2d+3} d\sqrt{2}\dst\max\left\{2^{2\n+6}d,\mathsf{C}_4 /2\right\}\mathsf{C}_5,\quad \mathsf{C}_8 \coloneqq 3\cdot 2^{3\n+1}\n^\n\ex^{-\n} d\mathsf{C}_4\sqrt{2} ,\\
& \mathsf{C}_9  \coloneqq 3\cdot 2^{-(4\n+2d)} d\mathsf{C}_4\sqrt{2}  +2^{-\n}\n^\n\ex^{-\n}\mathsf{C}_7\mathsf{C}_8,\quad  \mathsf{C}_{10}  \coloneqq 2\left(\frac{3}{2}\right)^{\n+d+1}\dst\int_{\rn} |y|_1^{\n+1}\ex^{-|y|_1}dy,\\
& \mathsf{C}_{11}  \coloneqq 8\left(\frac{3}{2}\right)^{3\t+d+2}\dst\int_{\rn}( 2|y|_1^{\t}+3|y|_1^{2\t+1}+|y|_1^{3\t+2})\ex^{-|y|_1}dy\;,\quad  \mathsf{C}_{12}  \coloneqq \max\{2 \mathsf{C}_{10} \;,2 \mathsf{C}_{11} \;,12\mathsf{C}_{0}\}\,,
\end{align*}
\begin{align*}
&\mathsf{C}_* \coloneqq \max\{ 2^{11\n+6d+4}\n^\n\ex^{-\n}\mathsf{C}_5^2\mathsf{C}_6\mathsf{C}_7\mathsf{C}_8,\; (2^{\n/2-2d+2}(18d^3+70)\n^\n\ex^{-\n}\mathsf{C}_{4} \mathsf{C}_9  \mathsf{C}_{12} ^{-1})^2,\; 2^{\t+8}d^2\mathsf{C}_{4} \mathsf{C}_9 \sqrt{2}\}\,.
\end{align*}
Then,
\beqa{constants}
&&c_*\coloneqq \mathsf{C}_*\,,\phantom{AAAAAA.}
c_0\coloneqq 2^{-4} \mathsf{C}_4\,,\quad\quad\phantom{.}
c_1\coloneqq  2^{\t-1/2}d^{-1} \mathsf{C}_9 \,,  
\nonumber
\\
&& 
c_2\coloneqq 2^{2\n+6}d \mathsf{C}_9  \,, \qquad 
c_3\coloneqq \mathsf{C}_9\,, 
\phantom{AAAAA.}
c_4\coloneqq \ex^{(4d)^{-1}}\mathsf{C}_{4} \mathsf{C}_9  \mathsf{C}_{12}^{-1}\,.
\eeqa

\erem

\rem{ArnPntwis}
There is a small flaw in \cite{chierchia2019ArnoldKAM}: The parameter\footnote{In the present remark, we will adopt the notations of \cite{chierchia2019ArnoldKAM}.} $\mathsf{L}$ chosen in \cite[Lemma~1]{chierchia2019ArnoldKAM} is not big enough to ensure that the new perturbation $P'$ and the symplectic change of coordinates $\phi$ are well--defined on $D_{\bar{r}/2,s'}(\mathscr D_\sharp')$.  The right choice is the following
\begin{align*}
& \mathsf{L}\coloneqq \mathsf{P}\dst\max\Big\{\frac{40d\mathsf{T}^2\mathsf{K}  }{r\bar{r}\s^{\n+d}}\,,\,\frac{2\mathsf{C}_4}{\a\bar{r}\s^{2(\n+d)}}\Big\}\;,\quad \mathsf{W}\coloneqq \diag(\bar{r}^{-1}\uno_d,\;\uno_d),\quad \hat{\epsilon}_0\coloneqq \mathsf{C}_{9}\;\s_0^{-2(\n+d)-1}\epsilon_0\;\torsion_0^2\;\l_0^\n\;,\\
& \mathsf{L}_j\coloneqq \frac{\mathsf{P}_j}{r_{j+1}}\dst\max\Big\{\frac{80d\sqrt{2}\mathsf{T}_0\th_0  }{r_j \s_j^{\n+d}}\,,\,\frac{\mathsf{C}_4}{\a \s_j^{2(\n+d)}}\Big\}\;,\quad \mathsf{W}_j\coloneqq \diag(2r_{j+1}^{-1}\uno_d,\;\uno_d),\quad \hat{\epsilon}_{j+1}\coloneqq \frac{\mathsf{K}_0\vae^{2^{j+1}}\mathsf{P}_{j+1}}{\a^2}\;,\\
& \mathsf{P}_{j+2}\coloneqq \l_*\th_*^{j+1}\frac{\mathsf{K}_0\mathsf{P}_{j+1}^2}{\a^2}\;,\quad \hat{\epsilon}_{j+1}\coloneqq \l_*\th_*^{j+2}\epsilon_{j+1}\;.
\end{align*} 
Of course, one needs, then, to change accordingly (and in a straightforward way) the constants involved, as follows:

\beqano
\n			 &\coloneqq& \t+1\;\\
\mathsf{C}_0 &\coloneqq& 4\sqrt{2}\left(\frac{3}{2}\right)^{2\n+d}\dst\int_{\rn} \left( |y|_1^{\n}+|y|_1^{2\n}\right)\ex^{-|y|_1}dy\;,\quad
\mathsf{C}_1 \coloneqq 2\left(\frac{3}{2}\right)^{\n+d}\dst\int_{\rn} |y|_1^{\n}\ex^{-|y|_1}dy\;,\\
\mathsf{C}_2 &\coloneqq& 2^{3d}d\;,\quad
\mathsf{C}_3 \coloneqq	\left(d^2\mathsf{C}_1^2+6d\mathsf{C}_1 +\mathsf{C}_2\right)\sqrt{2}\;,\quad
\mathsf{C}_4 \coloneqq \max\left\{6d^2\mathsf{C}_0,\,\mathsf{C}_3\right\}\;,\\
\mathsf{C}_5   &\coloneqq& \frac{3\cdot 2^5d}{5}\;,\quad
\mathsf{C}_6 \coloneqq \dst{\max}\left\{2^{2\n}\,,\,\mathsf{C}_5\right\}\;,\quad
\mathsf{C}_7 \coloneqq 3d\cdot 2^{4\n+2}\sqrt{2}\dst\max\left\{640d^2\,,\,\mathsf{C}_4 \right\}\;,\\
\mathsf{C}_8 &\coloneqq& \left(2^{-d}\mathsf{C}_6\right)^{1/8}\;,\quad
\mathsf{C}_{9} \coloneqq 3d\cdot 2^{2\n+2}\sqrt{2}\dst\max\left\{80d\sqrt{2}\,,\,\mathsf{C}_4\right\}\;,\\
\mathsf{C}_{10} &\coloneqq& (4\n\ex^{-1})^{2\n}\left(1+2^{4\n+2d+2}(\n\ex^{-1})^{2\n}\mathsf{C}_6^2\mathsf{C}_7 \right)\mathsf{C}_9/(3d^2)\;,\quad \mathsf{C}_{11} \coloneqq ({5d\cdot 2^{3(\n+1)}})^{-1}{\mathsf{C}_{10}}\;,\\
\mathsf{C}_{12} &\coloneqq& 2^{2(5\n+4d+2)}\;\mathsf{C}_6^2\;\mathsf{C}_7\;\mathsf{C}_8\;\mathsf{C}_9\;,\quad \mathsf{C}_{13} \coloneqq  \mathsf{C}_{10}+\mathsf{C}_{11}\;,\quad
\mathsf{C}_{14}\coloneqq  \mathsf{C}_{12}\;,\\
\mathsf{C}_{15}&\coloneqq& 18d^3+70\;, \qquad\mathsf{C}_{16}\coloneqq (6\n\ex^{-1})^{4\n}\;,\quad \mathsf{C}\ \coloneqq \max\{3\mathsf{C}_{10},\; \mathsf{C}_{13}\}\;,\\
\mathsf{C}_*&\coloneqq&  \max\left\{\mathsf{C}_{16}\mathsf{C}_{14}^{2/3},\;  6\mathsf{C}_{15}\mathsf{C}_{16}\mathsf{C}^2,\;2^{2(4\n+2d+1)}\mathsf{C}_{16}\mathsf{C}_9^2,\; \mathsf{C}_{10}^2%,\;(\n\ex^{-1})^{\n}\mathsf{C}_{11}
\right\}\;.
\eeqano
The smallness condition $(14)$ and the estimate $(16)$ become, respectively,
$$
\a\le \frac{r}{\mathsf{T}}\qquad \mbox{and}\qquad \epsilon\le \epsilon_*\coloneqq \frac{(s-s_*)^{a}}{\mathsf{C}_*\;\theta^6}\;,
$$
and
$$
\max\Big\{ \| u_*\|_{s_*}\,,\  \|\partial_x u_*\|_{s_*}\,,\,\frac{\mathsf{K}}{\a}\;(\log\epsilon^{-1})^\n\, \|v_*\|_{s_*}
\Big\} \le \frac{\mathsf{C}\  \torsion^3}{(s-s_*)^{a/2}} \ \epsilon\;(\log\epsilon^{-1})^\n \le \frac{1}{4\ex}\,,
$$
where $a\coloneqq 6\n+3d+2$.
\erem
\appB{Tools}
{\subsection{Classical estimates (Cauchy, Fourier)} 
\lemtwo{Cau} {\rm \cite{CC95}} 
Let $p\in \natural,\,r,s>0, y_0\in \cn$ and $f$ a real--analytic function $\cball_{r,s}(y_0)$ with 
$
\|f\|_{r,s}\coloneqq \sup_{\cball_{r,s}(y_0)}|f|<\infty.
$ 
Then,\\
{\bf (i)} For any multi--index $(l,k)\in \natural^d\times\natural^d$ with $|l|_1+|k|_1\le p$ and for any $0<r'<r,\, 0<s'<s$,\footnote{As usual, $\dpr_y^l\coloneqq \frac{\dpr^{|l|_1}}{\dpr y_1^{l_1}\cdots\dpr y_d^{l_d}},\, \forall\, y\in\rn,\, l\in\zn $.\label{notDevPart}}
\[\|\partial_{y}^l \partial_{x}^k f\|_{r',s'}\leq p!\; \|f\|_{r,s}(r-r')^{|l|_1}(s-s')^{|k|_1}.\]
{\bf (ii)}  For any $ k\in \zn$ and any $y\in \cball_r(y_0)$ 
$$|f_k(y)|\leq \ex^{-|k|_1 s}\|f\|_{r,s}.
$$
\elem
}
\subsection{An Inverse Function Theorem}

\thm{IFT}
Let  $D$ be a convex subset of $\cn$, $y_0\in D$  and let $f\in C^1(D,\cn)$ such that\footnote{$f'$ being the Jacobian matrix of $f$.}  $\det f'(y_0)\not=0$.Assume
\beq{IFT1}
 \varrho\coloneqq \sup_{y\in D} \|\uno-Tf'(y)\|<1\,,\qquad \qquad T\coloneqq (f'(y_0))^{-1}\,.
\eeq
Then, $\det f'(y)\not=0$, for each $y\in D$ and
\beq{InvJa}
\|(f'(y))^{-1}\|\le \l\coloneqq \frac{\|T\|}{1-\varrho}\,.
\eeq
Moreover, $f$ in injective on $D$ and its inverse function $g:f(D)\stackrel{\tiny \rm onto}{\to} D$ satisfies
\beq{lipg}
 \Lip_{f(D)}(g)\le \l\,.
 \eeq
Furthermore, if $D\coloneqq {B_r(y_0)}$, $\r\coloneqq r/\l$ and $z_0\coloneqq f(y_0)$, then 
\beq{gontof}
{B_\r(z_0)}\subset f(D)\,.
\eeq
\ethm

\proof 
For every $y\in D$, we have
$
f'(y)=f'(y_0)(\uno-A), 
$
where $A\coloneqq \uno- Tf'(y)$ with $\|A\|\le \varrho<1$. Thus, $f'(y)$ is invertible and
$$
\|(f'(y))^{-1}\|=\|(\sum_{n\ge 0} A^n)T\|\le \frac{\|T\|}{1-\varrho},
$$ 
proving \equ{InvJa}. Now, consider the auxiliary map $F\colon D\ni y\longmapsto y-Tf(y)$. We have $F\in C^1(D,\cn)$ and 
$
\sup_D \|F'\|\leby{IFT1} \varrho.
$
Thus, for every $y,\by\in D$ with $y\neq \by$, we have 
\beqa{IFT4}
\|T\| \|f(y)-f(\by)\|& \stackrel{\equ{IFT1}}{\ge} & \big\|T\big(f(y)-f(\by)\big)\big\|\nonumber\\
&=& \|(y-\by)+(F(\by)-F(y))\|\nonumber\\
&\ge&\|y-\by\|- \|y-\by\|\sup_D \|F'\|\nonumber\\
&\ge& \|y-\by\|(1-\varrho)\gtby{IFT1}0\,,
\eeqa
which shows that $f$ is injective on $D$ and, hence, that \equ{lipg} holds.
\\
To show \equ{gontof} in the case $D\coloneqq {B_r(y_0)}$ and
 $\r\coloneqq r/\l$, fix $\eta\in \cn$ with $\|\eta-z_0\|< \r$. We have to show that there exists $\by\in D$ such that $f(\by)=\eta$.
Define the map
\beq{Phi}
\Phi:y\in D \mapsto \Phi(y)\coloneqq y - T \big(f(y)-\eta\big)\in Y\,.
\eeq
Then, $\Phi$ is a contraction on $D$. Indeed, $\Phi$ is $C^1$, $\Phi'(y)=\uno-Tf'(y)$ and 
\beqa{Pcontr}
\Lip_D \Phi=\sup_D \|\Phi'\|=\varrho<1.
\eeqa
Furthermore, $\Phi:D\to D$, since, if $y\in D$, then
\beqano
\|\Phi(y)-y_0\|&\le& \|\Phi(y)-\Phi(y_0)\|+\|\Phi(y_0)-y_0\|\\
&\stackrel{\equ{Pcontr}}{\le}& \varrho\, r+ \|T\|\|\eta - z_0\|< \varrho\, r+\|T\|\, \r=r\,.
\eeqano
Hence, by the contraction Lemma, $\Phi$ has a (unique) fixed point $\by\in D$, but $\Phi(\by)=\by$ means $f(\by)=\eta$.
\qed

\subsection{Internal coverings}\label{app:covering}

\noi
Given any non--empty subset $D$ of $\rn$, and given $r>0$, a {\bf $r$--internal covering of $D$} is a subset $P$ of $D$
such that $D\subset \bigcup_{y\in P} \rball_r(y)$; the {\bf $r$--internal covering number of $D$}, denoted $N_r^{\rm int}(D)$,
is the minimal cardinality of any $r$--internal cover.

\Giu
In 
\cite{biasco2018explicit} the following simple upper bound (having fixed the sup norm in $\rn$) on $N_r^{\rm int}(D)$ for bounded sets $D$ is given:

\begin{lemma}\label{covering}
Let $D\subseteq\real^d$  be a non--empty bounded set. Then, for any $r>0$, one has\footnote{$[x]$ denotes the integer--part (or ``floor'') function $\max\{n\in \integer |\, \ n\le x\}$, while
$\lceil x\rceil$ denote the ``ceiling function''
$\min\{n\in \integer |\, \ n\ge x\}$; observe that $\lceil x\rceil\le [x]+1$.
}
\beq{ppp}
N_r^{\rm int}(D) \le \Big(\Big[ \frac{\diam D}{r}\Big]+1\Big)^d\,.
\eeq
\end{lemma}

\nl
For convenience of the reader, we reproduce here the elementary proof of the lemma.

\proof It is enough to produce a $r$--internal cover of $D$ with cardinality $N$ bounded by the right hand side of \equ{ppp}.\\
If $D$ is a singleton, the claim is obvious with $N=1$. Assume, now, $\d\coloneqq \diam D>0$,
and let $M\coloneqq [\d/r]+1$ and
 $z_i=\inf\{x_i|\ x\in D\}$. Then, $D\subseteq K:= z+[0,\delta]^d$ and one can find  $0<r'<r$ close enough to $r$ so that $\lceil\d/r'\rceil\le [\d/r]+1=M$. Then, one can 
cover $K$ with $M^d$ closed, contiguous cubes $K_j$, $1\le j\le M^d$, with edge of length $r'$. Let $j_i$ be the indices such that $K_{j_i}\cap D\neq \emptyset$ and pick a $y_i\in K_{j_i}\cap E$;
let $1\le N\le M^d$ be the number of such cubes.
Observe that, since we have chosen the sup--norm in $\real^d$, one has $K_{j_i}\subseteq \rball_{r}(y_i)$ and 
 \equ{ppp} follows. \qed

%%%%%%%%%%%%
\subsection{Extensions of Lipschitz continuous functions}
Here we recall a Theorem due to Minty according to which a Lipschitz continuous functions can be extended keeping unchanged both the sup--norm {\sl and} the Lipschitz constant.
\thmtwo{Minty}{G.~J. Minty\cite{minty1970extension}}
Let $(V,\average{\cdot\;,\cdot})$ be a separable inner product space, 
$\emptyset\neq A\subseteq V$, $L>0,\;0<\a\le 1$ and $g\colon A\to \rn$ a $(L,\a)$--Lipschitz--H\"older continuous function on $A$, namely, $g$ satisfies
\beq{aLipHol}
|g(x_1)-g(x_2)|_2\le L\; \|x_1-x_2\|^\a\;, \qquad\forall\; x_1,x_2\in A\;,
\eeq
where $\|\cdot\|$ denotes the norm on $V$ induced by the inner product.
Then, there exists a global $(L,\a)$--Lipschitz--H\"older continuous function\footnote{I.e.,  satisfying \equ{aLipHol} on $V$.} $G\colon V\to \rn$ such that $G|_A=g$. Futhermore, $G$ can be chosen in such away that $G(V)$ is contained in the closed convex hull of $g(A)$. Hence, in particular,
\beq{MinTextMint}
\dst\sup_{x\in V} |G(x)|_2=\dst\sup_{x\in A} |g(x)|_2\quad \mbox{and}\quad 
\dst\sup_{x_1\neq x_2\in V} \frac{|G(x_1)-G(x_2)|_2}{\|x_1-x_2\|^\a}=\dst\sup_{x_1\neq x_2\in A} \frac{|g(x_1)-g(x_2)|_2}{\|x_1-x_2\|^\a}\;.
\eeq
\ethm

\subsection{Lebesgue measure and Lipschitz continuous map \label{appC}}
\lem{LebLipLem}
Let $\emptyset\not=A\subset\rn$ be a Lebesgue--measurable set and $f\colon A\to \rn$ be Lipschitz continuous.
Then,
\beq{upperbound}
\meas\big(f(A)\big)\le \Lip_A(f)^d  \meas(A)
\eeq
and\footnote{Inequality \equ{EstApd} is sharp as shown by the example $f=(1+\d)\;\id$.} 
\beq{EstApd}
|\meas(f(A))-\meas(A)|\le ((1+\d)^d-1)\meas(A)\;.
\eeq
where
\beq{apceq1}
\d\coloneqq \Lip_A(f-\id)
\eeq
\elem
\proof
Eq. \equ{upperbound} is standard: see, e.g., Theorem 2, Sec 2.2 and Theorem 1, Sec 2.4 in \cite{evans2015measure}.

\nl
Let us prove \equ{apceq1}.
By Theorem~\ref{Minty}, $f-\id$ can be extended to a Lipschitz continuous $g\colon \rn\righttoleftarrow$ with
$$
\Lip(g) =\Lip_A(f-\id)= \d\;. 
$$
By Rademacher's Theorem, there exists a set $N\subset\rn$ with $\meas(N)=0$  such that $g$ is differentiable on $\rn\setminus N$ and  
$$
\|g_y\|_{\rn\setminus N}\le\Lip_{\rn\setminus N}(g)\le \Lip(g)=\d\;. %\le \su2\;.
$$
Now pick $y\in \rn\setminus N$. Then, 
\begin{align*}
|\det(\uno_d+g_y(y))-1|&=\left|\dst\int_0^1\frac{d}{dt}\det(\uno_d+tg_y)dt \right|=\left|\dst\int_0^1\tr\left(\adj(\uno_d+tg_y)g_y\right) dt \right|\\
					&\le \dst\int_0^1d\|\uno_d+tg_y\|^{d-1}\|g_y\|dt\le \dst\int_0^1 d\left(1+\d t\right)^{d-1}\d dt= (1+\d)^d-1.
\end{align*}
Thus, by the change of variable (or area) formula\footnote{See \cite{evans2015measure}, $\S3.3$.}, we have
\begin{align*}
|\meas(f(A))-\meas(A)|&=|\meas({(\id + g)}  (A))-\meas(A)|=\left|\dst\int_{(\id+g)(A)}dy-\int_A dy \right|\\
					  &= \left|\dst\int_{(\id+g)(A\setminus N)}dy-\int_{A\setminus N} dy \right|= \left|\dst\int_{A\setminus N}|\det (\uno_d+g_y)| dy-\int_{A\setminus N} dy \right|\\
					  &\le \dst\int_{A\setminus N}|\det(\uno_d+ g_y)-1| dy\le ((1+\d)^d-1)\meas(A)\;. \qedeq
\end{align*}

\subsection{Lipeomorphisms ``close'' to identity\label{appG}}

\lem{LipRang}
Let $g\colon \cn\to\cn$ be a Lipschitz continuous function such  that
\begin{align}
&\d\coloneqq \sup_{\rn}|g-\id|<\infty\;,\label{EqApG1}\\
&\th\coloneqq\Lip_\rn(g-\id)<1.\label{EqApG2}
\end{align}
Then, $g$ has a Lipschitz global inverse $G$ satisfying
\begin{align}
&\sup_{\rn}|G-\id|\le \d\;,\label{GG}\\
&\Lip_\rn(G-\id)\|<\frac1{1-\th}\label{GGG}\,.
\end{align}
Furthermore,  for any 
$\emptyset\neq A\subset\cn$,
%$A$, non--empty, closed subset of $\rn$,
\beq{gonto}
A\subset g\left(\ovl {\cball_{\d}( A)}\right)\;.
\eeq
%\eeq
\elem
\proof Let $f\coloneqq g-\id$, then, for any $x_i\in\rn$, one has
\beqano |g(x_1)-g(x_2)|&=& \big|x_1-x_2 + \big(f(x_1)-f(x_2)\big)\big|
\stackrel{\equ{EqApG2}}{>}
|x_1-x_2|-\th |x_1-x_2|\\
&=&(1-\th) |x_1-x_2|\,, 
\eeqano
which proves injectivity of $g$ and that
\beq{ginj}
\inf_{x_1\neq x_2} \frac{|g(x_1)-g(x_2)|}{|x_1-x_2|}\ge 1-\th>0\,.
\eeq
Let us now prove \equ{gonto}.
Let $\bar y\in A$. It is enough to show that there exists $|y|\le \d$ such that $\bar y=g(y+\bar y)$ \ie $y=-f(y+\bar{y})$ \ie $y$ is a fixed point of the map 
$$h\colon \ovl{\cball_{\d}(0)}\ni y\mapsto -f(y+\bar{y}).$$
But, for any $y\in \ovl{\cball_{\d}(0)}$,
$$
|h(y)|=|f(y+\bar{y})|\le \|f\|_{\rn}\leby{EqApG1} \d\;,
$$
\ie $h\colon \ovl{\cball_{\d}(0)}\to \ovl{\cball_{\d}(0)}$. Moreover, $h$ is a contraction since $\Lip_{\ovl{\cball_{\d}(0)}}(h)\le 
\Lip_\rn (f)\ltby{EqApG2}1$. Thus, by Banach's fixed point Theorem, we see that \equ{gonto} holds.
\\
From \equ{gonto} it follows at once that $g$ is onto $\rn$. \\
Now, \equ{GG} and \equ{GGG} follows easily from, respectively  \equ{EqApG1} and \equ{ginj}. 
\qed

\subsection{Whitney smoothness\label{appD}}
\dfn{WhitDef}
Let $A\subset \rn$ be non--empty and $n\in\natural_0$, $m\in\natural$. A function $f\colon A\to \real^m$ is said $C^n$ on $A$ in the Whitney sense, with Whitney derivatives $(f_\n)_{\n\in \natural_0^d,{|\n|_1\le n}}$ , $f_0=f$, and we write $f\in C^n_W(A,\real^m)$, if for any $\vae>0$ and $y_0\in A$, there exists $\d>0$ such that, for any $y,y'\in A\cap \rball_\d(y_0)$ and $\n\in \natural_0^d$, with  ${|\n|_1\le n}$,
\beq{whitdef}
\Big|f_\n(y')-\dst\sum_{\substack{\m\in \natural_0^d\\ {|\m|_1\le n-|\n|_1}}} \frac{1}{\m!}f_{\n+\m}(y)(y'-y)^\m\Big|\le \vae |y'-y|^{n-|\n|_1}\;.
\eeq
\edfn
\lemtwo{Whit1}{\cite{chierchia1986quasi,koudjinan2019quantitative}}
Let $A\subset \rn$ be non--empty and $n\in\natural_0$. For $m\in\natural$, let $f_m$ be  a real--analytic function with holomorphic extension to $D_{r_m}(A)$, with $r_m\downarrow0$ as $m\rightarrow\infty$. Assume that
\beq{hypWhit}
a\coloneqq\dst\sum_{m=1}^\infty \|f_m\|_{r_m,A}\;r_m^{-n}<\infty,\qquad \|f_m\|_{r_m,A}\coloneqq \dst\sup_{\rball^d_{r_m}(A)}|f_m|\;.
\eeq
Then $f\coloneqq\dst\sum_{m=1}^\infty f_m\in C^n_W(A,\real)$ with Whitney derivatives $f_\n\coloneqq\dst\sum_{m=1}^\infty\dpr_y^\n f_m$.
\elem
For completeness, we recall the beautiful Whitney extension theorem.
\thmtwo{Whi34}{\cite{whitney1934analytic}}
Let $A\subseteq \rn$ be a closed set and $f\in C^n_W(A,\real)$, $n\in\natural_0$. Then there exists $\bar{f}\in C^n(\rn,\real)$, real--analytic on $\rn\setminus A$ and such that $D^\n \bar{f}=f_\n$ on $A$, for any $\n\in \natural_0^d$, with ${|\n|_1\le n}$.
\ethm

\subsection{Measure of tubular neighbourhoods of hypersurfaces\label{appE}}

\nl
Recall  the definitions of minimal focal distance and of inner domains given in \S~\ref{sec:sm}.

\nl
The first elementary remark is that, for smooth domains,  taking $\r$--inner domains is the inverse operation of taking  $\r$--neighborhood:

\lem{lem:geo1}Let $\actiondom\subset \rn$ be an open and  bounded set with $C^2$ boundary $\partial \actiondom=S$  compact and connected. Then, for any $0<\r'<\r\le \minfoc(S)$, one has
\beq{minfoc1}
{\bf B}_\r\big(\actiondom''_\r\big)=\actiondom \,,\qquad\mbox{and}\qquad {\bf B}_{\r-\r'}\big(\actiondom''_\r\big)= \actiondom''_{\r'} \,.
\eeq
\elem

\proof
We start proving the first part of \equ{minfoc1}. By definition, ${\bf B}_\r\big(\actiondom''_\r\big)\subseteq\actiondom$. Thus, it remains only to show that $\actiondom\setminus \actiondom''_\r\subseteq {\bf B}_\r\big(\actiondom''_\r\big)$.

\begin{comment}
 Observe that for any given $y_0\in \dpr\actiondom''_\r$,
\beq{distBdryDr}
\dist_2(y_0,\rn\setminus \actiondom\big)=\dist_2(y_0,S\big)=\r.
\eeq
Indeed,  by definition, for any $y\in \rn\setminus \actiondom$, we have $y\not\in \rball_\r(y_0)$ \ie $|y-y_0|_2\ge \r$ \ie $\r'\coloneqq \dist_2(y_0,\rn\setminus \actiondom\big)\ge \r$. By contradiction, assume that $\r'> \r$. Let $y_1\in \ovl{\rball_{(\r'-\r)/2}(y_0)}$. We have, for any $y\in \rball_{\r}(y_1)$, 
$$
\dist_2(y,\rn\setminus \actiondom\big)\ge\dist_2(y_0,\rn\setminus \actiondom\big)- |y_0-y_1|_2-|y_1-y|_2\ge\r'-\frac{\r'-\r}{2}-\r=\frac{\r'-\r}{2}>0,
$$
and, as $\rn\setminus \actiondom$ is closed, then, $y\not\in\rn\setminus \actiondom$ \ie $y\in \actiondom$. Hence, $\rball_{\r}(y_1)\subseteq \actiondom$ which, by definition, means $y_1\in \actiondom''_\r$. Thus, $\ovl{\rball_{(\r'-\r)/2}(y_0)}\subseteq \actiondom''_\r$ and, in particular, $y_0$ is an interior point of $\actiondom''_\r$, contradiction. Thus, \equ{distBdryDr} holds.
\end{comment}
\noi
Let then $y_0\in \actiondom\setminus \actiondom''_\r$. As $S$ is compact and $\dist_2$ is continuous, there exists $\bar{y}_0\in S$ such that $\dist_2(y_0,\rn\setminus \actiondom\big)=\dist_2(y_0,S\big)=|y_0-\bar{y}_0|_2$. The vector $\n\coloneqq (y_0-\bar{y}_0)/|y_0-\bar{y}_0|_2$ is the inward unit normal to $\dpr\actiondom=S$ at $\bar{y}_0$. Indeed, for any smooth curve $\g\colon [0,1]\to S$ with $\g(0)=\bar{y}_0$, $0$ is a minimum of the smooth map $f(t)\coloneqq |\g(t)-y_0|_2^2$. Thus,
$$
0=f'(0)=2\dot{\g}(0)\cdot (\bar{y}_0-y_0).
$$
which, by the arbitrariness of $\g$, implies that the line $(\bar{y}_0y_0)$ is perpendicular to the tangent space to $S$ at $\bar{y}_0$ and, therefore $\n$ is the inward unit normal to $\dpr\actiondom$ at $\bar{y}_0$. Let $y_1\coloneqq \bar{y}_0+\r\n$. By assumption, we have $\dist_2(y_1,S\big)=\r$, and, therefore, $y_1\in \actiondom$. In addition, $y_1\in \actiondom''_\r$. Indeed, for any $y\in{\bf B}_\r(y_1)$, $\dist_2(y,\rn\setminus \actiondom\big)\ge \dist_2(y_1,\rn\setminus \actiondom\big)-|y_1-y|_2=\dist_2(y_1,S\big)-|y_1-y|_2=\r-|y_1-y|_2>0$. Thus, as $\rn\setminus \actiondom$ is a closed set,  $y\not\in \rn\setminus \actiondom$ \ie $y\in\actiondom$. Hence, ${\bf B}_\r(y_1)\subseteq \actiondom$ \ie $y_1\in \actiondom''_\r$. In particular, the argument above shows that:\footnote{Actually, %$y_1\in \dpr\actiondom''_\r$. Indeed, 
one checks easily that $\dpr\actiondom''_\r=\{y\in\rn\;:\; \dist_2(y,\rn\setminus\actiondom\big)=\r\}$ and $\inter(\actiondom''_\r)=\{y\in\rn\;:\; \dist_2(y,\rn\setminus \actiondom\big)>\r\}$, $\inter(\actiondom''_\r)$ being the interior of $\actiondom''_\r$.} for any $y\in \rn$, $\dist_2(y,\rn\setminus \actiondom\big)\ge \r$ implies that $y\in \actiondom''_\r$. Thus, as $y_0\in \actiondom\setminus \actiondom''_\r$, we have $\dist_2(y_0,\rn\setminus \actiondom\big)< \r$, which means $y_0$ is in the open segment $(\bar{y}_0,y_1)$. Therefore, $|y_0-y_1|_2< |\bar y_0-y_1|_2=\r$ \ie $y_0\in{\bf B}_\r(y_1)\subseteq {\bf B}_\r(\actiondom''_\r)$.

\noi
We now prove the second part of \equ{minfoc1}. We have ${\bf B}_{\r-\r'}\big(\actiondom''_\r\big)\subseteq \actiondom''_{\r'}$. Indeed, for any $y_0\in \actiondom''_{\r}$, $y_1\in {\bf B}_{\r-\r'}(y_0)$ and $y\in {\bf B}_{\r'}(y_1)$,
$$
|y-y_0|\le |y-y_1|+|y_1-y_0|< \r'+(\r-\r')=\r\quad \ie \ y\in {\bf B}_{\r}(y_0),
$$
which implies ${\bf B}_{\r-\r'}\big(\actiondom''_\r\big)\subseteq \actiondom''_{\r'}$. It remains to show that $\actiondom''_{\r'}\bks\actiondom''_{\r}\subseteq {\bf B}_{\r-\r'}\big(\actiondom''_\r\big)$. The proof follows in analogous to the previous. Let $y_0\in \actiondom''_{\r'}\bks\actiondom''_{\r}$ and $\bar{y}_0\in S$ such that $\dist_2(y_0,\rn\setminus \actiondom\big)=\dist_2(y_0,S\big)=|y_0-\bar{y}_0|_2$. Then, $\r'\le |y_0-\bar{y}_0|_2< \r$, and the  vector $\n\coloneqq (y_0-\bar{y}_0)/|y_0-\bar{y}_0|_2$ is the inward unit normal to $\dpr\actiondom=S$ at $\bar{y}_0$. Set $y_1'\coloneqq \bar{y}_0+\r'\n$. Thus, $|y_1'-\bar{y}_0|_2=\r'\le |y_0-\bar{y}_0|_2$ and, hence, $y_1'\in \actiondom''_{\r'}$ and $y_1'$ is in the semi open segment $(\bar{y}_0,y_0]$. Therefore,
$|y_1'-{y}_0|_2=|y_0-\bar{y}_0|_2-|y_1'-\bar{y}_0|_2<\r-\r'$. Hence, $y_0\in {\bf B}_{\r-\r'}(y_1')\subseteq {\bf B}_{\r-\r'}(\actiondom''_{\r'})$ \ie $\actiondom''_{\r'}\bks\actiondom''_{\r}\subseteq {\bf B}_{\r-\r'}\big(\actiondom''_\r\big)$.
\qed
\nl
Next result, gives a precise evaluation of tubular domains {\sl in the case the metric is the euclidean one}. Define
\beq{def.t.n}
{ \mathfrak T}_{\r}(S)\coloneqq  \{u\in \rn:\ \dist_2(u,S)< \r\}\;.
\eeq

\lem{lem:geo2} Let $\actiondom\subset \rn$ be a bounded set with $C^2$ boundary $\partial \actiondom=S$  compact and connected. Then, for any $0<\r\le \minfoc(S)$, then, 
\beqa{minfoc2}
\meas ({\mathfrak T}_{\r}(S))\le {\frac2d}\, \frac{(1+\r\kappa)^{d}-1}{\kappa}\;\mathcal{H}^{d-1}(S)\,,
\eeqa
where $\kappa\coloneqq \sup_S\max_{ 1\le j\le d-1}|\kappa_j|$ with $\kappa_j$ the principal curvatures of $S$, while $\mathcal{H}^{d-1}$ denotes the $(d-1)$--dimensional Hausdorff measure (`surface area').
\elem
\proof\!\footnote{Compare \cite{sternberg2012curvature}, Ch. 1.} 
We will estimate the `inner tubular neighbourhoods' 
$${\mathfrak T}_\r'(S)\coloneqq\{y\in \actiondom: \dist_2(y,S)<\r\}\,,
$$ 
as the argument for `outer tubular neighbourhood' $\{y\notin \actiondom: \dist_2(y,S)<\r\}$, 
is completely analogous.\\
Since $S$ is compact and connected, we may 
assume that $S=f^{-1}(\{0\})$ with $f\in C^2(\real^d,\real)$ and $0$  a regular value for $f$. Set
 $$
 \nu(x)=\frac{\nabla f}{|\nabla f|_2}\;,\qquad |\cdot|_2\coloneqq \dist_2(\cdot, 0)\;,
 $$ 
  and replacing eventually $f$ by $-f$, we can assume that $\nu$ is the inwards unit normal vector fields of $S$. Let  $\{\phi_j\colon U_j\to \real^m\}_{j=1}^p$ be an atlas of $S$,
 $$
 \Psi_j(u,t)\coloneqq \phi_j(u)+t \nu(\phi_j(u)),\qquad O_j\coloneqq \Psi_j(U_j\times [0,\r)),
 $$ and observe that\footnote{As $S=\bigcup_{j=1}^p \phi_j(U_j)$, we have $\mathscr T_{\r}(S)=\bigcup_{j=1}^p O_j$, for any $0<\r\le \minfoc(S)$.}
 $$
 {\mathfrak T}_{\r}'(S)=\bigcup_{j=1}^p O_j.
 $$ 
Let $\{\psi_j\}_{j=1}^p$ be a partition of unity subordinated to the open covering of $\{{O}_{j}\}_{j=1}^p$ of ${\mathfrak T}'_{\r}(S)$ \ie
\begin{itemize}
\item[$\mathbf{(i)}$] $\psi_j\in C^\infty_c( {\mathfrak T}'_{\r}(S))$ ;
\item[$\mathbf{(ii)}$] $0\le \psi_j\le 1$ ;
\item[$\mathbf{(iii)}$] $\supp\psi_j\subset {O}_{j}$ ;
\item[$\mathbf{(iv)}$] $\dst\sum_{j=1}^p\psi_j\equiv 1$ on  $ {\mathfrak T}'_{\r}(S)$\;.
\end{itemize}
Given $1\le j\le p$, define $n_j\colon  U_j\longrightarrow \mathbb{S}^d=\{x\in \real^d\;:\; |x|_2=x_1^2+\cdots+x_d^2=1\}\subset \real^d$ as
 $$
n_j\coloneqq   \nu\circ \phi_j,
 $$
and $K_j\colon U_j\longrightarrow T^*S$ such that\footnote{$T^*S$ being the cotangent bundle of $S$.},
 $$
 K_j(u)\coloneqq -\nu'(\phi_j(u)).
 $$
 Then, $K_j$ is symmetric\footnote{$K_j$ is actually the Weingarten map $\mathcal{W}_x=-\nu'(x)$ ``written in the local chart'' $(U_j,\phi_j)$.} and therefore diagonalizable, with eigenvalues 
 $\kappa_i\circ \phi_j^{-1}$, $1\le i\le d-1$  and satisfies
 \beq{Weiga}
 {\frac{\dpr n_j}{\dpr  u}}=-K_j \frac{\dpr \phi_j}{\dpr  u}\;.
 \eeq
Thus, recalling that $0=\dpr_x \nu^2=2\nu'\cdot \nu$, we have
 \begin{align*}
 \meas ({\mathfrak T}'_\r(S))&= \sum_{j=1}^p \int_{{O}_{j}}\psi_j\;du dt\\
                &= \sum_{j=1}^p \int_{\Psi_j(U_j\times [0,\r))}\psi_j\;du dt\\
                &= \sum_{j=1}^p \int_{U_j\times [0,\r)}\Psi_j^*(\psi_jdu dt)\\
                &= \sum_{j=1}^p \int_{U_j\times [0,\d)}\psi_j\circ \Psi_j\;\left|\det \left(\frac{\dpr \Psi_j}{\dpr(u,t)}\right)\right| du dt\\
                &\eqby{Weiga} \sum_{j=1}^p \int_{U_j\times [0,\r)}\psi_j\circ \Psi_j\;\left|\det \left[\frac{\dpr\phi_j}{\dpr u}-tK_j\cdot \frac{\dpr\phi_j}{\dpr u},  \;\nu(\phi_j(u))\right]\right| du dt
\\                &= \sum_{j=1}^p \int_{U_j\times [0,\r)}\psi_j\circ \Psi_j\;\left|\det \left(\begin{pmatrix}                
                \uno_{d-1}- t K_j               
                \end{pmatrix}\left[\frac{\dpr\phi_j}{\dpr u},  \;\nu(\phi_j(u))\right]\right)\right| du dt\\
                &= \sum_{j=1}^p \int_{U_j\times [0,\r)}\psi_j\circ \Psi_j\;\left|\det (  \uno_{d-1}- t K_j)\right|\left|\det \left[\frac{\dpr\phi_j}{\dpr u},  \;\nu(\phi_j(u))\right]\right| du dt\\
                &\le \int_0^\r \sum_{j=1}^p\int_{U_j}\psi_j\big(\phi_j(u)+t \nu(\phi_j(u))\big)\left|\det \left[\frac{\dpr\phi_j}{\dpr u},  \;\nu(\phi_j(u))\right]\right| du\; (1+t\kappa)^{d-1}  dt\\
                &{=\int_0^\r \sum_{j=1}^p\int_{U_j}\psi_j\big(\phi_j(u)+t \nu(\phi_j(u))\big)\left(\det \left(\frac{\dpr\phi_j}{\dpr u}\right)^T\frac{\dpr\phi_j}{\dpr u}\right)^{1/2}\; du\; (1+t\kappa)^{d-1}  dt}
                \end{align*}
\begin{align*}
\phantom{\meas ({\mathfrak T}'_\r(S))}                
                &{= \int_0^\r\sum_{j=1}^p \int_{\phi_j(U_j)}\psi_j\big(x+t \nu(x)\big)\;d\mathcal{H}^{d-1}(x)\; (1+t\kappa)^{d-1} dt} \quad \mbox{\footnotesize (see \cite[Theorem 2, pg. 99]{evans2015measure})}\\
                &\stackrel{\mathbf{(ii)}}{\le} \int_0^\r\sum_{j=1}^p \int_{\bigcup_{i=1}^p\phi_i(U_i)}\psi_j\big(x+t \nu(x)\big)\;d\mathcal{H}^{d-1}(x)\; (1+t\kappa)^{d-1} dt\\
                &= \int_0^\r \int_{S}\sum_{j=1}^p\psi_j\big(x+t \nu(x)\big)\;d\mathcal{H}^{d-1}(x)\; (1+t\kappa)^{d-1} dt\\
                &\stackrel{\mathbf{(iv)}}{=}\int_0^\r \int_{S}d\mathcal{H}^{d-1}(x)\; (1+t\kappa)^{d-1} dt\\
                &= \frac{(1+\r\kappa)^{d}-1}{d\;\kappa}\;\mathcal{H}^{d-1}(S)\;. \qedeq
 \end{align*}
\subsection{Kolmogorov non--degenerate normal forms\label{app9}}
Let $H\colon \mathcal{M}\coloneqq \rn\times\tn\to\real$ be a $C^2$--Hamiltonian. An embedded torus $\mathcal{T}$ in $\mathcal{M}$ is said $H$--Kolmogorov non--degenerate if there exists a neighborhood $\mathcal{M}_0$ of $\{0\}\times\tn$ in $\mathcal{M}$, a symplectic change of coordinates $\phi\colon \mathcal{M}_0\to \mathcal{M}$ with $\phi(\{0\}\times\tn)=\mathcal{T}$, a constant $E\in \real$, a vector $\o\in\rn$ and a function $Q\colon \mathcal{M}_0\to\real$ of class $C^2$ such that
\beq{KolNomFrm}
H\circ \phi(y,x)=E+\o\cdot y+Q(y,x)\qquad\mbox{and}\qquad \dpr_y^\mu Q(0,\cdot)\equiv 0,\ \forall\; \mu\in\natural_0^d, \  |\mu|_1\le 1,
\eeq
and
\beq{KolmNrmF2}
\det\average{\dpr_{yy}Q(0,\cdot)}\not=0.
\eeq
A Hamiltonian $H$ in the form \equ{KolNomFrm} is said in Kolmogorov normal form. The Kolmogorov normal form is said non--degenerate if, in addition, the quadratic (in $y$) part $Q$ satisfies \equ{KolmNrmF2}.
\bibliographystyle{apa}
\bibliography{BibtexDatabase}

\begin{thebibliography}{Whi34}

\bibitem[Arg20]{Argentieri2020}
F.~Argentieri.
\newblock Isolated points of {D}iophantine sets, Preprint, 2020.

\bibitem[Arn63]{ARV63}
V.I. Arnold.
\newblock {Proof} of {A.N.} {Kolmogorov}'s theorem on the conservation of
  conditionally periodic motions with a small variation in the {Hamiltonian}.
\newblock {\em Russian Math. Surv}, 18(9), 1963.

\bibitem[BC15]{Biasco2015}
Luca Biasco and Luigi Chierchia.
\newblock On the measure of {L}agrangian invariant tori in nearly-integrable
  mechanical systems.
\newblock {\em Atti Accad. Naz. Lincei Rend. Lincei Mat. Appl.},
  26(4):423--432, 2015.

\bibitem[BC17]{biasco2017secondary}
Luca Biasco and Luigi Chierchia.
\newblock Kam 2017 theory for secondary tori.
\newblock {\em {\rm arXiv} (arXiv:1702.06480v1)}, 2017.

\bibitem[BC18]{biasco2018explicit}
L.~Biasco and L.~Chierchia.
\newblock {Explicit} estimates on the measure of primary {KAM} tori.
\newblock {\em Annali di Matematica Pura ed Applicata (1923-)},
  197(1):261--281, 2018.

\bibitem[BC20]{biasco2020nonlinearity}
L.~Biasco and L.~Chierchia.
\newblock On the topology of nearly-integrable {H}amiltonians at simple
  resonances.
\newblock {\em Nonlinearity}, 33(7):3526--3567, 2020.

\bibitem[CC95]{CC95}
A.~Celletti and L.~Chierchia.
\newblock A constructive theory of {Lagrangian} tori and computer-assisted
  applications.
\newblock In {\em Dynamics reported}, pages 60--129. Springer, 1995.

\bibitem[CG82]{chierchia1982smooth}
L.~Chierchia and G.~Gallavotti.
\newblock Smooth prime integrals for quasi-integrable {H}amiltonian systems.
\newblock {\em Nuovo Cimento B (11)}, 67(2):277--295, 1982.

\bibitem[Chi86]{chierchia1986quasi}
L.~Chierchia.
\newblock Quasi-periodic {Schroedinger} operators in one dimension, absolutely
  continuous spectra, bloch waves, and integrable hamiltonian systems.
\newblock Technical report, New York Univ., NY (USA), 1986.

\bibitem[CK19]{chierchia2019ArnoldKAM}
L.~Chierchia and C.E. Koudjinan.
\newblock {V.I.~Arnold's} ``pointwise'' {KAM Theorem}.
\newblock {\em Regular and Chaotic Dynamics}, 24(6):583--606, 2019.

\bibitem[EG15]{evans2015measure}
L.~C. Evans and R.~F. Gariepy.
\newblock {\em {Measure} {Theory} and fine properties of functions}.
\newblock CRC press, 2015.

\bibitem[Kou19]{koudjinan2019quantitative}
C.E. Koudjinan.
\newblock {\em Quantitative KAM normal forms and sharp measure estimates}.
\newblock PhD thesis, Universit\`a degli Studi Roma Tre, March 2019.
\newblock arxiv.org/abs/1904.13062.

\bibitem[Min70]{minty1970extension}
G.~J. Minty.
\newblock On the extension of {Lipschitz}, {Lipschitz}--{H{\"o}lder}
  continuous, and monotone functions.
\newblock {\em Bulletin of the American Mathematical Society}, 76(2):334--339,
  1970.

\bibitem[Nei81]{neishtadt1981estimates}
A.I. Neishtadt.
\newblock Estimates in the {Kolmogorov} theorem on conservation of
  conditionally periodic motions.
\newblock {\em Journal of Applied Mathematics and Mechanics}, 45(6):766--772,
  1981.

\bibitem[P{\"o}s82]{poschel1982integrability}
J~P{\"o}schel.
\newblock Integrability of hamiltonian systems on cantor sets.
\newblock {\em Communications on Pure and Applied Mathematics}, 35(5):653--696,
  1982.

\bibitem[P{\"o}s01]{JP}
J.~P{\"o}schel.
\newblock A {Lecture} on the {Classical} {KAM} {Theorem}.
\newblock In {\em Proc. Symp. Pure Math}, volume~69, pages 707--732, 2001.

\bibitem[Sal04]{salamon2004kolmogorov}
D.~Salamon.
\newblock The {Kolmogorov}--{Arnold}--{Moser} {Theorem}.
\newblock {\em Math. Phys. Electron. J}, 10(3):1--37, 2004.

\bibitem[Ste12]{sternberg2012curvature}
S.~Sternberg.
\newblock {\em Curvature in Mathematics and Physics}.
\newblock Courier Corporation, 2012.

\bibitem[Whi34]{whitney1934analytic}
H.~Whitney.
\newblock Analytic extensions of differentiable functions defined in closed
  sets.
\newblock {\em Transactions of the American Mathematical Society},
  36(1):63--89, 1934.

\bibitem[Zeh10]{zehnder2010lectures}
E.~Zehnder.
\newblock {\em Lectures on dynamical systems: {Hamiltonian} vector fields and
  symplectic capacities}, volume~11.
\newblock European Mathematical Society, 2010.

\end{thebibliography}

\end{document}